\documentclass{amsproc}

\usepackage{amsmath}
\usepackage{amsfonts}
\usepackage{amsthm}

\newcommand{\baton}[1]{\mathbb #1}
\newcommand{\E}{{\baton E}}
\newcommand{\N}{{\baton N}}
\newcommand{\R}{{\baton R}}
\newcommand{\Z}{{\baton Z}}
\newcommand{\Q}{{\baton Q}}
\newcommand{\T}{{\baton T}}

\newcommand{\CG}{{\mathcal G}}
\newcommand{\CI}{{\mathcal I}}
\newcommand{\CJ}{{\mathcal J}}

\newcommand{\CX}{{\mathcal X}}
\newcommand{\CY}{{\mathcal Y}}
\newcommand{\CZ}{{\mathcal Z}}
\newcommand{\CF}{{\mathcal F}}

\newcommand{\CH}{{\mathcal H}}

\newcommand{\CB}{{\mathcal B}}
\newcommand{\CC}{{\mathcal C}}

\newcommand{\norm}[1]{\lVert #1\rVert}
\newcommand{\nnorm}[1]{\lvert\!|\!| #1|\!|\!\rvert}
\newcommand{\inv}{^{-1}}
\newcommand{\wh}{\widehat}
\newcommand{\diffsym}{\raise -1mm\hbox{$\Delta$}}
\newcommand{\one}{\boldsymbol{1}}

\newcommand{\type}[1]{^{[#1]}}
\newcommand{\bn}{\mathbf{n}}
\newcommand{\bzero}{\mathbf{0}}
\newcommand{\bx}{\mathbf{x}}
\newcommand{\by}{\mathbf{y}}
\newcommand{\beps}{\mathbf{\epsilon}}
\newcommand{\tf}{\tilde{f}}
\newcommand{\limproj}{\underleftarrow\lim\,}

\newtheorem{theorem}{Theorem}[section]
\newtheorem{lemma}[theorem]{Lemma}
\newtheorem{corollary}[theorem]{Corollary}
\newtheorem{proposition}[theorem]{Proposition}

\theoremstyle{definition}

\newtheorem{example}[theorem]{Example}

\theoremstyle{remark}

\numberwithin{equation}{section}

\begin{document}

\title{Ergodic methods in additive combinatorics}

\author{Bryna Kra}
\address{Department of Mathematics, Northwestern University,
2033 Sheridan Road,  Evanston, IL 60208-2730}
\email{kra@math.northwestern.edu}
\thanks{The first author was supported in part by NSF Grant DMS-\#0555250.}

\subjclass{Primary 37A30; Secondary 11B25, 27A45}

\keywords{Ergodic theory, additive combinatorics}

\begin{abstract}

Shortly after Szemer\'edi's proof that a set of 
positive upper density contains arbitrarily long 
arithmetic progressions, Furstenberg gave a new proof 
of this theorem using 
ergodic theory.  This gave rise to the field of 
ergodic Ramsey Theory, in which problems 
motivated by additive combinatorics are proven using 
ergodic theory.  Ergodic Ramsey Theory has since produced
combinatorial results, some of which have yet 
to be obtained by other means, and has also given 
a deeper understanding of the structure of 
measure preserving systems.  We outline the ergodic 
theory background needed to understand these results, 
with an emphasis on recent developments in 
ergodic theory and the relation to recent developments 
in additive combinatorics.  

These notes are based on four lectures given during 
the School on Additive Combinatorics at the 
Centre de Recherches Math\'ematiques, Montreal in 
April, 2006.  The talks were aimed at an audience without 
background in ergodic theory.  
No attempt is made to include complete proofs of all statements 
and often the reader is referred to the original sources.  Many of 
the proofs included are classic, included 
as an indication of which ingredients play a role in the 
developments of the past ten years.
\end{abstract}

\maketitle

\section{Combinatorics to ergodic theory}
\subsection{Szemer\'edi's Theorem}
Answering a long standing conjecture of
Erd\H os and Tur\'an~\cite{ET}, Szemer\'edi~\cite{S} 
showed that a set $E\subset\Z$ with positive upper 
density\footnote{Given a set $E\subset\Z$,  its 
{\em upper density} $d^*(E)$ is defined by 
$d^*(E) = \limsup_{N\to\infty}\frac{|E\cap\{1, \ldots, N\}|}{N}$.}
contains arbitrarily long arithmetic progressions.  Soon thereafter, 
Furstenberg~\cite{F} gave a new proof of 
Szemer\'edi's Theorem using ergodic theory, 
and this has lead to the rich field of ergodic 
Ramsey theory.  Before describing some of the results in this subject, 
we motivate the use of ergodic theory for 
studying combinatorial problems.  

We start with the finite formulation of Szemer\'edi's Theorem:
\begin{theorem}[Szemer\'edi~\cite{S}]
Given $\delta > 0$ and $k\in\N$, there is a function 
$N(\delta, k)$ such that if $N > N(\delta, k)$ 
and $E\subset\{1, \ldots, N\}$ is a subset 
with $|E|\geq \delta N$, then $E$ 
contains an arithmetic progression of length $k$.  
\end{theorem}

It is clear that this statement immediately implies the first
formulation 
of Szemer\'edi's Theorem, and a compactness argument gives 
ths converse implication.  

\subsection{Translation to a probability system}
Starting with Szemer\'edi's Theorem, one gains insight 
into the intersection 
of sets in a probability system\footnote{By a {\em probability 
system}, we mean a triple $(X, \CX, \mu)$ where $X$ 
is a measure space, $\CX$ is a $\sigma$-algebra of measurable 
subsets of $X$, and $\mu$ is a probability measure.  
In general, we use 
the convention of denoting the $\sigma$-algebra $\CX$ 
by the associated calligraphic 
version of the measure space $X$.}:

\begin{corollary}  
\label{cor:probsystem}
Let $\delta > 0$, $k\in\N$, 
$(X, \CX, \mu)$ be a probability space 
and $A_1, \ldots, A_N\in\CX$  with $\mu(A_i)\geq\delta$
for $i=1, \ldots, N$.  If $N> N(\delta, k)$, 
then there exist $a, d\in \N$ such that 
$$
A_a\cap A_{a+d}\cap A_{a+2d}\cap\ldots\cap A_{a+kd}\neq\emptyset \ .
$$
\end{corollary}

\begin{proof}
For $A\in\CX$, let 
$\one_A(x)$ denote the characteristic function of $A$ (meaning 
that $\one_A(x)$ is $1$ for $x\in A$ and is $0$ otherwise). Then 
$$
\int_X \frac{1}{N}\sum_{n=0}^{N-1}\one_{A_n}\, d\mu \geq \delta \ .
$$

Thus there exists $x\in X$ with 
$\frac{1}{N}\sum_{n=0}^{N-1}\one_{A_n}(x) \geq \delta$.  
Then $E = \{n\colon x\in A_n\}$
satisfies $|E|\geq \delta N$, and so Szemer\'edi's Theorem 
implies that $E$ contains an arithmetic progression of length $k$.  
\end{proof}

\subsection{Measure preserving systems}

A {\em probability measure preserving system} is a 
quadruple $(X, \CX, \mu, T)$, where $(X, \CX, \mu)$ is a 
probability space and $T\colon X\to X$ is 
a bijective, measurable, measure preserving transformation.  
This means that for all $A\in\CX$, $T^{-1}A\in\CX$ and 
$$\mu(T^{-1}A) = \mu(A)\ . $$
In general, we refer to a probability measure 
preserving system as a {\em system}. 

Without loss of generality, we can 
place several simplifying assumptions on our systems.  
We assume that $\CX$ is countably generated; thus  for 
$1 \leq p < \infty$, $L^p(\mu)$ is separable.  
We implicitly assume that all sets and functions are measurable
with respect to the appropriate $\sigma$-algebra, even 
when this is not explicitly stated. 
Equality between sets or functions is meant up to sets 
of measure $0$.  

\subsection{Furstenberg multiple recurrence}
In a system, one can use Szemer\'edi's Theorem to derive a 
bit more information about intersections of sets.  
If $(X, \CX, \mu, T)$ is a system and $A\in\CX$ with 
$\mu(A)\geq\delta > 0$, then 
$$A, T^{-1}A, T^{-2}A, \ldots, T^{-n}A, \ldots$$ 
are all sets of measure $\geq \delta$.  
Applying Corollary~\ref{cor:probsystem} to this sequence of sets,
we have the existence of $a, d\in \N$ with
$$
T^{-a}A\cap T^{-(a+d)}A\cap T^{-(a+2d)}\cap\ldots \cap
T^{-(a+kd)}A\neq\emptyset \ .
$$

Furthermore, the measure of this intersection must be positive.  If not, 
we could remove from $A$ a subset of measure zero containing all the 
intersections and obtain a subset of measure at least $\delta$ 
without this property.  In this way, 
starting with Szemer\'edi's Theorem, 
we have derived Furstenberg's multiple recurrence theorem:

\begin{theorem}[Furstenberg~\cite{F}]
\label{th:multrec}
Let $(X,\CX,\mu,T)$ be a system, and let $A\in\CX$ with $\mu(A)>0$.  
Then for any $k\geq 1$, there exists $n\in\N$ such that 
\begin{equation}
\label{eq:multrec}
\mu\bigl(A\cap T^{-n}A\cap T^{-2n}A
\cap\dots\cap 
 T^{-kn}A\bigr) > 0 \ .
\end{equation}
\end{theorem}

\section{Ergodic theory to combinatorics}
\subsection{Strong form of multiple recurrence}
We have seen that 
Furstenberg multiple recurrence can be easily derived from 
Szemer\'edi's Theorem.  More interesting is the converse 
implication, showing that one can use ergodic theory to prove 
regularity properties of subsets of the integers. 
This approach has two major components, 
and has been since used to deduce other patterns 
in subsets of integers with positive upper density.  
(See Section~\ref{sec:other}.)  
The first is proving a certain recurrence statement in ergodic theory,
like that of Theorem~\ref{th:multrec}.  The second is 
showing that this statement implies a corresponding statement 
about subsets of the integers.  We now make this more precise.

To use ergodic theory to show that 
some intersection of sets has positive measure, 
it is natural to average the expression 
under consideration.  This leads us to the strong form 
of Furstenberg's multiple recurrence:

\begin{theorem}[Furstenberg~\cite{F}]
\label{th:strongmultrec}
Let $(X,\CX,\mu,T)$ be a system and let $A\in\CX$ with $\mu(A)>0$.
Then for any $k\geq 1$,  
\begin{equation}
\label{eq:strongmultrec}
\liminf_{N\to\infty}\frac{1}{N}\sum_{n=0}^{N-1}\mu(A\cap T^{-n}A\cap 
T^{-2n}A\cap\ldots\cap T^{-kn}A)
\end{equation}
is positive.  
\end{theorem}

In particular, this implies the existence of infinitely many 
$n\in\N$ such that the intersection in~\eqref{eq:multrec} is 
positive and Theorem~\ref{th:multrec} follows.  
We return later to a discussion of how to prove 
Theorem~\ref{th:strongmultrec}.  

\subsection{The correspondence principle}
The second major component is using this multiple 
recurrence statement to derive a statement about 
integers, such as Szemer\'edi's Theorem.  
This is the content of Furstenberg's Correspondence 
Principle:

\begin{theorem}[Furstenberg~\cite{F},~\cite{Fbook}]
Let $E\subset\Z$ have positive upper density. There
exist a system $(X,\CX,\mu,T)$ and a set $A\in\CX$ with
$\mu(A)= d^*(E)$ such that
$$\mu( T^{-m_1}A\cap\dots\cap T^{-m_k}A)\leq
d^*\bigl( (E+m_1)\cap\dots\cap(E+m_k)\bigr)$$
for all $k\in\N$ and all $m_1,\dots,m_k\in\Z$.
\end{theorem}

\begin{proof}
Let  $X = \{0,1\}^\Z$ be endowed with the product topology and the
shift map $T$ given by $Tx(n)=x(n+1)$ for all $n\in\Z$.  
A point of $X$ is thus a sequence $x = \{x(n)\}_{n\in\Z}$, and the distance between 
two points $x = \{x(n)\}_{n\in\Z}, y = \{y(n)\}_{n\in\Z}\in X$ is defined to be $0$ if $x = y$ 
and $2^{-k}$ if $x\neq y$ and 
$k = \min\{|n|\colon x(n)\neq y(n)\}$.  
Define $a\in \{0,1\}^\Z$ by 
\begin{equation*}
a(n) = \begin{cases}\displaystyle
   1 &\text{ if }n\in E\\
   0&\text{ otherwise } \end{cases}       
\end{equation*}
and let $A=\{x\in X\colon x(0)=1\}$. Thus $A$ 
is a clopen (closed and open) set.  

For all $n\in\Z$, 
$$
T^{-n}a\in A\text{  if and only if }n\in E\ .
$$

By definition of $d^*(E)$, there exist 
sequences $\{M_i\}$ and $\{N_i\}$ of integers with 
$N_i\to\infty$ such that 
$$
\lim_{i\to\infty}\frac{1}{N_i}\bigl\vert E\cap [M_i, M_i+N_i-1]\bigr\vert
\to d^*(E)\ .
$$

Then
$$
\lim_{i\to\infty}\frac{1}{N_i}\sum_{n=M_i}^{M_i+N_i-1}\one_A(T^na) 
= \lim_{i\to\infty}\frac{1}{N_i}\sum_{n=M_i}^{M_i+N_i-1}\one_E(n) = 
d^*(E) \ .
$$

Let $\CC$ be the countable algebra generated by cylinder sets, 
meaning sets that are defined by specifying finitely many coordinates 
of each element and leaving the others free.
We can define an additive 
measure $\mu$ on $\CC$ by
$$
\mu(B) = \lim_{i\to\infty}\frac{1}{N_i}
\sum_{n=M_i}^{M_i+N_i-1}\one_B(T^na)\ ,
$$
where we pass to subsequences $\{N_i\}$, $\{M_i\}$  such that 
this limit exists for all $B\in\CC$.  (Note that $\CC$ is 
countable and so by diagonalization we can arrange it 
such that this limit exists for all elements of $\CC$.)

We can extend the additive measure to a $\sigma$-additive 
measure $\mu$ on all Borel sets $\CX$ in $X$, which is exactly 
the $\sigma$-algebra generated by $\CC$.  
Then $\mu$ is an invariant measure, meaning that for all $B\in\CC$,
$$
\mu(T^{-1}B) = 
\lim_{i\to\infty}\frac{1}{N_i}
\sum_{n=M_i}^{M_i+N_i-1}\one_B(T^{n-1}a) = \mu(B) \ .
$$
Furthermore, 
$$
\mu(A) = \lim_{i\to\infty}\frac{1}{N_i}
\sum_{n=M_i}^{M_i+N_i-1}\one_A(T^na) = d^*(E) \ .
$$
If $m_1, \ldots, m_k\in\Z$, then the set $T^{-m_1}A\cap \ldots\cap T^{-m_k}A$ 
is a clopen set, its indicator function is continuous, and 
\begin{multline*}
\mu(T^{-m_1}A\cap \ldots\cap T^{-m_k}A) 
 = \lim_{i\to\infty}\frac{1}{N_i}
\sum_{n=M_i}^{M_i+N_i-1}\one_{T^{-m_1}A\cap \ldots\cap T^{-m_k}A}(T^na) \\
 = \lim_{i\to\infty}\frac{1}{N_i}
\sum_{n=M_i}^{M_i+N_i-1}\one_{(E+m_1)\cap \ldots\cap (E+m_k)}(n)
\leq d^*\bigl((E+m_1)\cap \ldots\cap (E+m_k)\bigr)\ .
\end{multline*}
\end{proof}

We use this to deduce Szemer\'edi's Theorem from Theorem~\ref{th:multrec}.  
As in the proof of the Correspondence Principle, define 
$a\in \{0,1\}^\Z$ by 
\begin{equation*}
a(n) = \begin{cases}\displaystyle
   1 &\text{ if }n\in E\\
   0&\text{ otherwise } \ ,  \end{cases}       
\end{equation*}
and set $A=\{x\in \{0,1\}^\Z\colon x(0)=1\}$.  
Thus $T^na\in A$ if and only if $n\in E$.  

By Theorem~\ref{th:multrec}, there exists $n\in\N$ such that 
$$
\mu(A\cap T^{-n}A\cap T^{-2n}A\cap\ldots\cap T^{-kn}A) > 0\ .
$$
Therefore for some $m\in\N$, $T^m a$ enters this multiple intersection 
and so 
$$
a(m) = a(m+n) = a(m+2n) = \ldots = a(m+kn) = 1 \ . 
$$
But this means that 
$$
m, m+n, m+2n, \ldots, m+kn\in E
$$
and so we have an arithmetic progression of length $k+1$ 
in $E$.

\section{Convergence of multiple ergodic averages}
\subsection{Convergence along arithmetic progressions}

Furstenberg's multiple recurrence theorem left open the question 
of the existence of the limit in~\eqref{eq:strongmultrec}.  
More generally, one can ask if given a system 
$(X, \CX, \mu, T)$ and $f_1, f_2, \ldots, f_k\in L^\infty(\mu)$, does 
\begin{equation}
\label{eq:APav}
\lim_{N\to\infty}\frac{1}{N}\sum_{n=0}^{N-1}f_1(T^nx)\cdot f_2(T^{2n}x)\cdot 
\ldots \cdot f_k(T^{kn}x)
\end{equation}
exist?  
Moreover, we can ask in what sense (in $L^2(\mu)$ or pointwise) does this 
limit exist, and if it does exist, what can be said about the 
limit?  
Setting each function $f_i$ to be the indicator 
function of a measurable set $A$, we are back in the 
context of Furstenberg's Theorem.  

For $k=1$, existence of the limit in $L^2(\mu)$ 
is the mean ergodic theorem of von Neumann.  In Section~\ref{sec:vonneumann}, 
we give a proof of this statement.  For $k=2$, existence 
of the limit in $L^2(\mu)$ was proven by Furstenberg~\cite{F} as part of his 
proof of Szemer\'edi's Theorem.  
Furthermore, in the same paper he showed the existence of the limit 
in $L^2(\mu)$ in a weak mixing system for arbitrary $k$; we 
define weak mixing in Section~\ref{sec:wm} and outline the proof 
for this case.  

For $k\geq 3$, the proof requires a more subtle understanding 
of measure preserving systems, and we begin discussing this 
case in Section~\ref{sec:higher}.  Under some technical hypotheses, 
the existence of the limit in $L^2(\mu)$ for $k=3$ was first 
proven by Conze and Lesigne (see~\cite{CL3} and~\cite{CL4}), then by 
Furstenberg and Weiss~\cite{FW}, and in the general case by Host and 
Kra~\cite{HK1}.  More generally, we showed the existence of the 
limit for all $k\in\N$:
\begin{theorem}[Host and Kra~\cite{HK4}]
\label{th:mainAP}
Let $(X, \CX, \mu, T)$ be a system, let $k\in\N$,  and 
let $f_1, f_2, \ldots, f_k\in L^\infty(\mu)$.  Then 
the averages 
$$
\frac{1}{N}\sum_{n=0}^{N-1}f_1(T^nx)\cdot f_2(T^{2n}x)\cdot 
\ldots \cdot f_k(T^{kn}x)
$$
converge in $L^2(\mu)$ as $N\to\infty$.  
\end{theorem}

Such a convergence result for a finite system is trivial.  
For example, if $X = \Z/N\Z$, then $\CX$ consists of all partitions of $X$
and $\mu$ is the uniform probability measure, meaning that 
the measure of a set is proportional to the cardinality of a set.  
The transformation $T$ is given by $Tx = x+1\mod N$.  
It is then trivial to check the convergence of the 
average in~\eqref{eq:APav}.  
However, although the ergodic theory is trivial in this case, 
there are common themes to be explored, 
and throughout these notes, an effort 
is made to highlight the connection with 
recent advances in additive combinatorics (see~\cite{Kra2} 
for more on this connection).  Of particular interest 
is the role played by nilpotent groups, and homogeneous spaces of 
nilpotent groups, in the proof of the ergodic statement.  
Some of these connections are further discussed in the notes 
of Ben Green and Terry Tao.  

Much of the present notes is devoted to understanding the ingredients 
in the proof of Theorem~\ref{th:mainAP}, and the 
role of nilpotent groups in this proof.  
Other expository accounts of this 
proof can be found in~\cite{Host} and in~\cite{Kra1}.
$2$-step nilpotent groups first appeared in the work of Conze-Lesigne 
in their proof of convergence for $k=3$, and a $(k-1)$-step nilpotent 
group plays a similar role in convergence for the average in~\eqref{eq:APav}.  
Nilpotent 
groups also play some role in the combinatorial setup, 
and this has been recently verified by 
Green and Tao (see~\cite{GT2},~\cite{GT3}, and~\cite{GT4}) 
for progressions of length $4$ (which 
corresponds to the case $k=3$ in~\eqref{eq:APav}).  
For more on this connection, see the lecture notes of 
Ben Green in this volume.  

\subsection{Other results}
Using ergodic theory, other patterns have been shown to exist 
in sets of positive upper density and we discuss 
these results in Section~\ref{sec:other}.
We briefly summarize these results.
A striking example is the 
theorem of Bergelson and Leibman~\cite{BL} showing the 
existence of polynomial patterns in such sets.  Analogous to 
the linear average corresponding to arithmetic progressions, 
existence of the associated polynomial averages has been 
shown in~\cite{HK5} and~\cite{Leibman1}.  One can also 
average along `cubes'; existence of these averages and a 
corresponding combinatorial statement was shown in~\cite{HK4}.  
For commuting transformations, little is known and these partial 
results are summarized in Section~\ref{sec:commuting}.  
An explicit formula for the limit in~\eqref{eq:APav} was 
given by Ziegler~\cite{ziegler1}, who also has recently 
given a new proof~\cite{ziegler2} of Theorem~\ref{th:mainAP}.

\section{Single convergence (the case $k=1$)}

\subsection{Poincar\'e Recurrence}
The case $k=1$ in Furstenberg's multiple recurrence 
(Theorem~\ref{th:multrec}) is Poincar\'e Recurrence:

\begin{theorem}[Poincar\'e~\cite{poincare}]
If $(X,\CX,\mu,T)$ is a system  
and $A\in\CX$ with $\mu(A)>0$, 
then there exist infinitely many $n\in\N$ such that 
$\mu(A\cap T^{-n}A)> 0$.  
\end{theorem}

\begin{proof}
Let $F = \{x\in A\colon T^{-n}x\notin A \text{ for all } n\geq 1\}$.
Assume that $F\cap T^{-n}F = \emptyset $ for all $n\geq 1$.  This implies 
that for all integers $n\neq m$, 
$$
T^{-m}A\cap T^{-n}A = \emptyset \ .
$$
In particular, $F, T^{-1}F, T^{-2}F, \ldots$ are all pairwise disjoint sets 
and each set in this sequence has measure 
equal to $\mu(F)$.  
If $\mu(F) > 0$, then 
$$
\mu\bigl(\bigcup_{n\geq 0}T^{-n}F\bigr) = \sum_{n\geq 0}\mu(F) = \infty \ , 
$$
a contradiction of $\mu$ being a probability measure.  

Therefore $\mu(F) = 0$ and the statement is proven.  
\end{proof}

In fact the same proof shows a bit more: by a simple modification 
of the definition of $F$, we have that 
$\mu$-almost every $x\in A$ returns to $A$ infinitely often.

\subsection{The von Neumann Ergodic Theorem}
\label{sec:vonneumann}

Although the proof of Poincar\'e Recurrence is simple, 
unfortunately there seems to be no way to 
generalize it for multiple recurrence. 
Instead we prove a stronger statement, taking the 
average of the expression under consideration and showing that the 
$\liminf$ of this average is positive.  It is not any harder (for $k=1$ only!) to 
show that the limit of this average exists (and is positive).  This 
is the content of the von Neumann mean ergodic theorem.  
We first give the statement in a general Hilbert space:

\begin{theorem}[von Neumann~\cite{vonneumann}]
\label{th:vn}
If $U$ is an isometry of a Hilbert space $\CH$ and 
$P$ is the orthogonal projection onto the $U$-invariant 
subspace $\CI = \{f\in\CH\colon Uf = f\}$, then for all $f\in\CH$, 
$$
\lim_{N\to\infty}\sum_{n=0}^{N-1}U^nf = Pf \ .
$$
\end{theorem}

Thus the case $k=1$ in Theorem~\ref{th:mainAP} is an 
immediate corollary of Theorem~\ref{th:vn}.

\begin{proof}
If $f\in \CI$, then 
$$
\frac{1}{N}\sum_{n=0}^{N-1}U^nf = f
$$
for all $N\in\N$ and so obviously the average 
converges to $f$.  
On the other hand, if $f = g-Ug$ for some $g\in\CH$, then 
$$\sum_{n=0}^{N-1}U^nf = g-U^Ng$$ 
and so the average converges to $0$ as $N\to\infty$.  
Setting $\CJ = \{g-Ug\colon g\in\CH\}$ and taking 
$f_k\in \CJ$ and $f_k\to f\in\overline{\CJ}$, then 
\begin{multline*}
\Bigl\|\frac{1}{N}\sum_{n=0}^{N-1}U^nf\Bigr\| \leq 
\Bigl\|\frac{1}{N}\sum_{n=0}^{N-1}U^n(f-f_k)\Bigr\| + 
\Bigl\|\frac{1}{N}\sum_{n=0}^{N-1}U^n(f_k)\Bigr\| 
\\
\leq 
\Bigl\|\frac{1}{N}\sum_{n=0}^{N-1}U^n\Bigr\|\cdot\Bigl\|f-f_k\Bigr\| 
+ \Bigl\|\frac{1}{N}\sum_{n=0}^{N-1}U^n(f_k)\Bigr\| \ .
\end{multline*}
Thus for $f\in\overline{\CJ}$, 
the average $\frac{1}{N}\sum_{n=0}^{N-1}U^nf$ 
converges to $0$ as $N\to\infty$.  

We now show that an arbitrary $f\in\CH$ can be written as a 
combination of functions which exhibit these behaviors, meaning 
that any $f\in\CH$ can be written as $f = f_1+f_2$ 
for some $f_1\in\CI$ and $f_2\in\overline{\CJ}$.  
If $h\in \CJ^\perp$, then for all $g\in\CH$, 
$$
0 = \langle h, g-Ug\rangle = \langle h, g\rangle - 
\langle h, Ug\rangle \\ = \langle h,g\rangle -\langle U^*h, g\rangle 
= \langle h-U^*h, g\rangle 
$$
and so $h = U^*h$ and $h = Uh$.  
Conversely,  reversing the steps we have that 
if $h\in \CI$, then $h\in \CJ^\perp$. 

Since $\overline{\CJ}^\perp = \CJ^\perp$, we have 
$$
\CH = \CI\oplus\overline{\CJ}\ .
$$

Thus writing $f = f_1+f_2$ with 
$f_1\in\CI$ and $f_2\in\overline{\CJ}$, 
we have 
$$
\frac{1}{N}\sum_{n=0}^{N-1}U^nf = \frac{1}{N}\sum_{n=0}^{N-1}U^nf_1 + \frac{1}{N}\sum_{n=0}^{N-1}U^nf_2 \ .
$$ 
The first sum converges to the identity and the second sum to $0$.  
\end{proof}

Under a mild hypothesis on the system, we have an explicit 
formula for the limit.  
Let $(X, \CX, \mu, T)$ be a system.  
A subset $A\subset X$ is {\em invariant} if $T^{-1}A = A$.  The 
invariant sets form a sub-$\sigma$-algebra $\CI$ of $\CX$.  
The system $(X, \CX, \mu, T)$ is said to be {\em ergodic} if 
$\CI$ is trivial, meaning that 
every invariant set has measure $0$ or measure $1$.  

A measure preserving transformation $T\colon X\to X$ 
defines a linear operator $U_T\colon L^2(\mu)\to L^2(\mu)$ 
by 
$$
(U_Tf)(x) = f(Tx) \ .
$$
It is easy to check that the operator $U_T$ is a unitary operator 
(meaning its adjoint is equal to its inverse).  In a standard abuse of 
notation, we use the same letter to denote the operator and 
the transformation, writing $Tf(x) = f(Tx)$ instead of 
the more cumbersome $U_Tf(x) = f(Tx)$.  

We have: 
\begin{corollary} 
\label{cor:two}
If $(X, \CX, \mu, T)$ is a system and $f\in L^2(\mu)$, then 
$$
\frac{1}{N}\sum_{n=0}^{N-1}f(T^nx)
$$ 
converges in $L^2(\mu)$, as $N\to\infty$, to a 
$T$-invariant function $\tilde f$.  
If the system is ergodic, then the 
limit is the constant function $\int f\, d\mu$.  
\end{corollary}

Let $(X, \CX, \mu, T)$ be an ergodic system and let 
$A, B\in\CX$.  
Taking $f = \one_A$ in Corollary~\ref{cor:two} 
and integrating with respect to $\mu$ over a set $B$, 
we have:
$$
\lim_{N\to\infty}\frac{1}{N}\sum_{n=0}^{N-1}\int_B\one_A(T^nx)\,d\mu(x) 
 = \int_B\bigl(\int\one_A(y)\,d\mu(y)\bigr)\,d\mu(x) \ .
$$
This means that 
$$
\lim_{N\to\infty}\frac{1}{N}\sum_{n=0}^{N-1}
\mu(A\cap T^{-n}B) = \mu(A)\mu(B) \ .
$$

In fact, one can check that this condition is equivalent 
to ergodicity.

As already discussed, convergence in the case of the finite system $\Z/N\Z$ with 
the transformation of adding $1\mod N$, is trivial.  Furthermore
this system is ergodic. 
More generally, any permutation on $\Z/N\Z$ can be expressed as a product 
of disjoint cyclic permutations.  These permutations 
are the `indecomposable' invariant subsets of an arbitrary 
transformation 
on $\Z/N\Z$ and the restriction of the transformation to one 
of these subsets is ergodic.  

This idea of dividing a space into indecomposable components
generalizes: an arbitrary 
measure preserving system can be decomposed into, perhaps 
continuously many, 
indecomposable components, and these are exactly the ergodic ones.   
Using this {\em ergodic decomposition} 
(see, for example,~\cite{CFS}), instead of working with 
an arbitrary system, we reduce most of 
the recurrence and convergence questions 
we consider here to the same problem in an ergodic system.

\section{Double convergence (the case $k=2$)}

\subsection{A model for double convergence}
We now turn to the case of $k=2$ in Theorem~\ref{th:mainAP}, 
and study convergence of the double average
$$
\frac{1}{N}\sum_{n=0}^{N-1}f_1(T^nx)\cdot f_2(T^{2n}x)
$$ 
for bounded functions $f_1$ and $f_2$.  
Our goal is to explain how a 
simple class of systems, the rotations, suffice to understand 
convergence for the double average.  

First we explicitly define what is meant by a rotation.  
Let $G$ be a compact abelian group, with Borel $\sigma$-algebra 
$\CB$, Haar measure $m$, and fix 
some $\alpha\in G$.  Define $T\colon G\to G$ by 
$$
Tx = x+\alpha \ .
$$
The system $(G, \CB, m, T)$ is called a {\em group rotation}.  
It is ergodic if and only if $\Z\alpha$ is dense in $G$.  
For example, when $X$ is the circle $\T = \R/\Z$ and $\alpha\notin\Q$, 
the rotation by $\alpha$ is ergodic.

The double average is the simplest example of a {\em nonconventional} 
ergodic average: even for an ergodic system, the limit is not 
necessarily constant.  This sort of behavior does not 
occur for the single average of von Neumann's Theorem, where 
we have seen that the limit is constant in 
an ergodic system. 
Even for the simple example of an 
an ergodic rotation, the limit of the double average
is not constant:

\begin{example}
\label{ex:rotation}
Let $X = \T$, with Borel $\sigma$-algebra 
and  Haar measure, and let $T\colon X\to X$ be the 
rotation $Tx = x+\alpha\mod 1$.  Setting 
$f_1(x) = \exp(4\pi i x)$ and $f_2(x) = \exp(-2\pi i x)$, 
then  for all $n\in\N$, 
$$
f_1(T^nx)\cdot f_2(T^{2n}x) = \overline{f_2(x)} \ .
$$
In particular, this double average converges to a nonconstant function.
\end{example}

More generally, if $\alpha\notin\Q$ and $f_1, f_2\in L^\infty(\mu)$, 
the double average converges to 
$$
\int_\T f_1(x+t)\cdot f_2(x+2t)\, dt \ .
$$

We shall see that Fourier analysis suffices to understand 
this average.  By taking both functions to be the indicator 
function of a set with positive measure and integrating 
over this set, we then have that Fourier analysis 
suffices for the study of arithmetic progressions 
of length $3$, giving a proof of Roth's Theorem 
via ergodic theory.  Later we shall see that other 
more powerful methods are needed 
to understand the average along longer progressions.  
In a similar vein, rotations are the model 
for an ergodic average with $3$ terms, but are 
not sufficient for more terms.  We introduce 
some terminology to make these notions more precise.  

\subsection{Factors}

For the remainder of 
this section, we assume that $(X, \CX, \mu, T)$ is an ergodic system.  

A {\em factor} of a system $(X, \CX, \mu, T)$ 
can be defined in one of several equivalent ways.  
It is a $T$-invariant sub-$\sigma$-algebra $\CY$ of $\CX$.  
A second characterization is that a factor is 
a system $(Y,\CY,\nu,S)$ and a 
measurable map $\pi:X\to Y$, the {\em factor map}, such that 
$\mu\circ\pi^{-1}=\nu$ and $S\circ\pi=\pi\circ T$ for $\mu$-almost 
every $x\in X$.  A third characterization is that a factor is a 
$T$-invariant subalgebra $\CF$ of $L^\infty(\mu)$.  
One can check that the first two definitions 
agree by identifying $\CY$ with 
$\pi\inv(\CY)$, and that the first and third agree by 
identifying $\CF$ with $L^\infty(\CY)$.  
When any of these conditions holds, we say that $Y$, 
or the appropriate sub-$\sigma$-algebra, is a factor 
of $X$ and write $\pi\colon X\to Y$ for the factor map.   
We make use of a slight (and standard) abuse of notation, useing 
the same letter $T$ to denote both the transformation 
in the original system and the transformation in 
the factor system.  
If the factor map $\pi\colon X\to Y$ 
is also injective, we say that the 
two systems $(X, \CX, \mu, T)$ and $(Y, \CY, \nu, S)$ 
are {\em isomorphic}.  

For example, if $(X, \CX, \mu, T)$ and $(Y, \CY, \nu, S)$ are systems, 
then each is a factor of the product system 
$(X\times Y, \CX\times \CY, \mu\times\nu, 
T\times S)$ and the associated factor map is projection 
onto the appropriate coordinate.

A more interesting example can be given in the system
$X = \T\times \T$, with Borel $\sigma$-algebra and Haar measure, and 
transformation $T\colon X\to X$ given by 
$$T(x,y) = (x+\alpha, y+x)\ .$$
Then $\T$ with the rotation $x\mapsto x+\alpha$ is a factor of $X$.

\subsection{Conditional expectation}

If $\CY$ is a $T$-invariant sub-$\sigma$-algebra of $\CX$ 
and $f\in L^2(\mu)$, the {\em conditional expectation
$\E(f\mid\CY)$ of $f$ with respect to $\CY$} is the function on $Y$
defined by $\E(f\mid Y)\circ\pi = \E(f\mid\CY)$.  
It is characterized as
the $\CY$-measurable function on $X$ such that 
$$
\int_X f(x)\cdot g(\pi(x))\,d\mu(x) = \int_Y\E(f\mid \CY)(y)\cdot 
g(y)\, d\nu(y)
$$
for all $g\in L^\infty(\nu)$ and satisfies the identities 
$$
\int \E(f\mid\CY)\,d\mu = \int f\,d\mu  
$$
and 
$$
T\E(f\mid\CY) = \E(Tf\mid\CY) \ .
$$

As an example, take $X = \T\times \T$ endowed with 
the transformation $(x,y) \mapsto (x+\alpha, y+x)$.  
We have a factor $Z=\T$ endowed with the map 
$x\mapsto x+\alpha$.  Considering 
$f(x,y) = \exp(x) + \exp(y)$, we have 
$\E(f\mid \CZ) = \exp(x)$.  
The factor sub-$\sigma$-algebra $\CZ$ is the $\sigma$-algebra of 
sets that depend only on the $x$ coordinate.

\subsection{Characteristic factors}

For $f_1, \ldots, f_k\in L^\infty(\mu)$, 
we are interested in convergence in $L^2(\mu)$ of: 
\begin{equation}
\label{eq:AP}
\frac{1}{N}\sum_{n=0}^{N-1}T^nf_1\cdot T^{2n}f_2\cdot\ldots\cdot T^{kn}f_k\ .
\end{equation}
Instead of working with the whole system $(X, \CX, \mu, T)$, 
it turns out that it is easier to find some factor of 
the system that characterizes this average, meaning that 
if we have some means of understanding 
convergence of the average 
under consideration in a well chosen factor, then we can 
also understand convergence of the same average 
in the original system.  This motivates the following 
definition.

A factor $Y$ of $X$ is {\em characteristic} for the 
average~\eqref{eq:AP} if this average converges to $0$ 
when $\E(f_i\mid\CY) = 0$ for some $i\in\{1, \ldots, k\}$.  
This is equivalent to showing that 
the difference between~\eqref{eq:AP} and 
$$
\frac{1}{N}\sum_{n=0}^{N-1}T^n\E(f_1\mid\CY)\cdot 
T^{2n}\E(f_2\mid\CY)\cdot \ldots 
\cdot T^{kn}\E(f_k\mid\CY)\cdot  
$$
converges to $0$ in $L^2(\mu)$.  

By definition, the whole system is 
always a characteristic factor.  
Of course nothing is gained by using 
such a characteristic factor, and the notion 
only becomes useful when we can find a 
characteristic factor that has 
useful geometric and/or algebraic properties.  
A very short outline of the proof of convergence 
of the average~\eqref{eq:AP} is as follows: find a 
characteristic factor that has sufficient 
structure so as to allow one to prove convergence.  
We return to this idea later.

The definition of a characteristic factor can be 
extended for any other average under consideration, 
with the obvious changes: the limit remains 
unchanged when each function is replaced by its 
conditional expectation on this factor.  
This notion has been implicit in the literature 
since Furstenberg's proof of Szemer\'edi's Theorem, 
but the terminology was only introduced more recently 
in~\cite{FW}.

\subsection{Weak mixing systems}
\label{sec:wm}

The system $(X, \CX, \mu, T)$ is {\em weak mixing} if for all 
$A, B\in\CX$, 
$$
\lim_{N\to\infty}\frac{1}{N}\sum_{n=0}^{N-1}\bigl\vert 
\mu(T^{-n}A\cap B) -\mu(A)\mu(B)\bigr\vert = 0 \ .
$$

There are many equivalent formulations of this property, 
and we give a few (see, for example~\cite{CFS}):
\begin{proposition}
Let $(X,\CX, \mu, T)$  be a system.  The following are equivalent:
\begin{enumerate}
\item $(X, \CX, \mu, T)$ is weak mixing.  

\item There exists $J\subset\N$ of density zero 
such that for all $A,B\in\CX$
$$\mu(T^{-n}A\cap B)\to \mu(A)\mu(B) 
\text{ as } n\to\infty \text{ and }n\notin J\ .$$ 

\item For all $A, B, C\in\CX$ with 
$\mu(A)\mu(B)\mu(C)> 0$, there exists $n\in\N$ such that 
$$\mu(A\cap 
T^{-n}B)\mu(A\cap T^{-n}C) > 0\ .
$$

\item The system 
$(X\times X, \CX\times\CX, \mu\times\mu, T\times T)$ is ergodic.  
\end{enumerate}
\end{proposition}

Any system exhibiting rotational behavior (for example a 
rotation on a circle, or a system with a nontrivial 
circle rotation as a factor) is not weak mixing.  
We have already seen in Example~\ref{ex:rotation} 
that weak mixing, or lack thereof, has an effect 
on multiple averages.  We give a second example 
to highlight this effect:
\begin{example}
Suppose that 
$X = X_{1}\cup X_{2}\cup X_{3}$ with $T(X_{1}) = X_{2}$, 
$T(X_{2}) = X_{3}$ and $T(X_{3}) = X_{1}$, and that $T^{3}$ restricted 
to $X_i$, for $i=1,2,3$, is weak mixing.  
For the double average 
$$
\frac{1}{N}\sum_{n=0}^{N-1}f_1(T^{n}x)\cdot f_2(T^{2n}x)\ , 
$$
where $f_1, f_2 \in L^{\infty}(\mu)$, if 
$x\in X_{1}$, this average converges to 
$$
\frac{1}{3}\Bigl(\int_{X_{1}}f_1\,d\mu\int_{X_{1}}f_2\,d\mu + 
\int_{X_{2}}f_1 \,d\mu\int_{X_{3}}f_2\,d\mu + \int_{X_{3}}f_1\,d\mu
\int_{X_{2}}f_2\,d\mu \Bigr) \ .
$$
A similar expression with obvious changes holds
for $x\in X_{2}$ or $x\in X_{3}$.  
\end{example}

The main point is that (for the double average) the answer 
depends on the rotational behavior of the system.  This example 
lacks weak mixing and so has nontrivial rotation factor.  
We now formalize this notion.

\subsection{Kronecker factor}

The {\em Kronecker factor} $(Z_1, \CZ_1, m, T)$ of $(X, \CX, \mu, T)$ is the 
sub-$\sigma$-algebra of $\CX$ spanned by the eigenfunctions.  
A classical result is that the Kronecker factor 
can be given the structure of a group rotation:
\begin{theorem}[Halmos and von Neumann~\cite{HVN}]  
The Kronecker factor of a system is isomorphic to a 
system $(Z_1, \CZ_1, m, T)$, where $Z_1$ is a compact abelian group, 
$\CZ_1$ is its Borel $\sigma$-algebra, $m$ is the Haar measure, 
and $Tx = x+\alpha$ for some fixed $\alpha\in Z_1$.  
\end{theorem}

We use $\pi_1\colon X\to Z_1$ to 
denote the factor map from a system $(X, \CX, \mu, T)$ 
to the Kronecker factor $(Z_1, \CZ_1, m, T)$.
Then any eigenfunction of $X$ takes the form 
$$
f(x) = c\gamma(\pi_1(x)) \ ,
$$ 
where $c$ is a constant 
and $\gamma\in\widehat{Z_1}$ is a character of $Z_1$.  

We give two examples of Kronecker factors:

\begin{example}
If $X = \T\times\T$, $\alpha\in\T$, and $T\colon X\to X$ is the map
$$
T(x,y) = (x+\alpha, y+x)\ ,
$$
then the rotation $x\mapsto x+\alpha$ on $\T$ is the Kronecker factor 
of $X$.  It corresponds to the pure point spectrum.  (The spectrum 
in the orthogonal complement of the Kronecker factor 
is countable Lebesgue.)
\end{example}

\begin{example}
If $X = \T^3$, $\alpha\in\T$, and $T\colon X\to X$ is the map
$$
T(x,y,z) = (x+\alpha, y+x, z+y)\ ,
$$
then again the rotation $x\mapsto x+\alpha$ on $\T$ is the Kronecker factor 
of $X$.  This example has the same 
pure point spectrum as the first example, but the first example 
is a factor of the second example.  
\end{example}

The Kronecker factor can be used to 
give another characterization of weak mixing:
\begin{theorem}  [Koopman and von Neumann~\cite{KVN}] 
A system is not weak mixing if and only if it has a 
nontrivial factor which is a rotation on a compact abelian group.  
\end{theorem}
The largest of these factors is the Kronecker factor.

\subsection{Convergence for $k=2$}
\label{sec:k3}

If we take into account the rotational behavior in a system, meaning 
the Kronecker factor, then 
we can understand the limit of the double average
\begin{equation}
\label{eq:doubleav}
\frac{1}{N}\sum_{n=0}^{N-1}T^{n}f_1\cdot T^{2n}f_2 \ .
\end{equation}

An obvious constraint is that for $\mu$-almost every $x$, the triple 
$(x, T^nx, T^{2n}x)$ projects to an arithmetic progression in the 
Kronecker factor $\CZ_1$.  Furstenberg proved that 
this obvious restriction is the only restriction, 
showing that to prove convergence of double average, one 
can assume that the system is an ergodic rotation on a compact 
abelian group:

\begin{theorem}[Furstenberg~\cite{F}]
\label{th:F2}
If $(X, \CX, \mu, T)$ is an ergodic system, $(Z_1, \CZ_1,$ $m, T)$ 
is its Kronecker factor, and $f_1, f_2, \in 
L^{\infty}(\mu)$, then the limit 
$$
\Bigl\Vert \frac{1}{N}\sum_{n=0}^{N-1}
T^nf_1\cdot T^{2n}f_2 -\frac{1}{N}\sum_{n=0}^{N-1}
T^n\E(f_1\mid\CZ_1)\cdot T^{2n}\E(f_2\mid\CZ_1)
\Bigr\Vert_{L^2(\mu)}$$
tends to $0$ as $N\to\infty$.
\end{theorem}

In our terminology, this theorem can be quickly summarized: 
the Kronecker factor is characteristic for the double average.
To prove the theorem, we use a standard trick for 
averaging, which is an iterated use of a 
variation of the van der Corput Lemma on 
differences (see~\cite{KN} or~\cite{Berg3}):  

\begin{lemma}[van der Corput]
\label{lem:vdc}
Let $\{u_n\}$ be a sequence in a Hilbert space with 
$\Vert u_n\Vert \leq 1 $ for all $n\in\N$.  For $h\in\N$, set 
$$\gamma_h = \limsup_{N\to\infty}\Bigl\vert \frac{1}{N}\sum_{n=0}^{N-1}
\langle u_{n+h}, u_n\rangle \Bigr\vert \ .
$$
Then 
$$
\limsup_{N\to\infty}\Bigl\Vert\frac{1}{N}\sum_{n=0}^{N-1} u_n
\Bigr\Vert^2 \leq \limsup_{H\to\infty}\frac{1}{H}\sum_{h=0}^{H-1}
\gamma_h \ .
$$
\end{lemma}

\begin{proof}
Given $\epsilon > 0$ and $M\in\N$, for $N$ sufficiently 
large we have that 
$$
\Bigl\vert \frac{1}{N}\sum_{n=0}^{N-1}u_n - \frac{1}{N}\frac{1}{H}
\sum_{n=0}^{N-1}\sum_{h=0}^{H-1}u_{n+h} \Big\vert < \epsilon \ .
$$
By convexity, 
$$
\Bigl\Vert \frac{1}{N}\sum_{n=0}^{N-1}\frac{1}{H}\sum_{h=0}^{H-1}
u_{n+h}\Bigr\Vert^2
\leq 
\frac{1}{N}\sum_{n=0}^{N-1}\Bigl\Vert 
\frac{1}{H}\sum_{h=0}^{H-1}u_{n+h}\Bigr\Vert^2 = 
\frac{1}{N}\frac{1}{H^2}\sum_{n=0}^{N-1}
\sum_{h_1, h_2=0}^{H-1}\langle u_{n+h_1}, u_{n+h_2}\rangle 
$$
and this approaches
$$
\frac{1}{H^2}\sum_{h_1, h_2}^{H-1}\gamma_{h_1-h_2}
$$ 
as $N\to\infty$.  But the assumption implies that 
this approaches $0$ as $H\to\infty$.  
\end{proof}

We now use this in the proof of Furstenberg's Theorem:
\begin{proof}[of Theorem~\ref{th:F2}]
Without loss of generality, we assume that $\E(f\mid\CZ_1) = 0$ 
and we show that the double averages converges to $0$.  
Set $u_n = T^nf_1\cdot T^{2n}f_2$.  Then
\begin{multline*}
\langle u_n, u_{n+h}\rangle = 
\int T^nf_1\cdot T^{2n}f_2\cdot T^{n+h}\overline{f_1}
\cdot T^{2n+2h}\overline{f_2} \, d\mu \\
=\int (f_1\cdot T^h \overline{f_1})\cdot T^n(f_2\cdot T^{2h}
\overline{f_2})\, d\mu  \ .
\end{multline*}
By the Ergodic Theorem, 
$$
\lim_{N\to\infty}\frac{1}{N}\sum_{n=0}^{N-1}
\langle u_n, u_{n+h}\rangle
$$ 
exists and is equal to 
\begin{equation}
\label{eq:int}
\int f_1\cdot T^h\overline{f_1} \cdot {\mathbb P}
(f_2\cdot T^{2h}\overline{f_2})\, d\mu \ ,
\end{equation}
where $\mathbb P$ is projection onto the $T$-invariant functions of $L^2(\mu)$.  
Since $T$ is ergodic, $\mathbb P$ is projection onto the constant functions.  
But since $\E(f_1\mid\CZ_1) = 0$, $f$ is orthogonal to the constant 
functions and so the integral in~\eqref{eq:int} is $0$.  
The van der Corput Lemma immediately gives the result.
\end{proof}

Furstenberg used a similar argument combined with 
induction to show that in a weak mixing 
system, the average~\ref{eq:AP} converges to the 
product of the integrals in $L^2(\mu)$ for all $k\geq 1$.  

Finally, to show that a set of integers with positive 
upper density contains arithmetic progressions of length three 
(Roth's Theorem), by Furstenberg's Correspondence Principle it 
suffices to show double recurrence:
\begin{theorem}[Theorem~\ref{th:multrec} for $k=2$]
Let $(X, \CX, \mu, T)$ be an ergodic 
system, and let $A\in\CX$ with $\mu(A) > 0$. 
There exists $n\in\N$ with 
$$\mu(A\cap T^{-n}A\cap T^{-2n}A)> 0 \ .
$$
\end{theorem}

\begin{proof}
Let $f = \one_A$.  Then 
$$
\mu(A\cap T^{-n}A\cap T^{-2n}A) = \int f\cdot T^nf\cdot T^{2n}f\, d\mu \ .
$$
It suffices to show that 
$$
\lim_{N\to\infty}\frac{1}{N}\sum_{n=0}^{N-1}\int f\cdot T^nf\cdot 
T^{2n}f\, d\mu
$$ 
is positive.  

By Theorem~\ref{th:F2}, the limiting behavior 
of the double average $\frac{1}{N}\sum_{n=0}^{N-1}T^nf\cdot T^{2n}f$ 
is unchanged if $f$ is replaced by $\E(f\mid\CZ_1)$.  
Multiplying by $f$ and integrating, it thus suffices to show that  
\begin{equation}
\label{eq:intint}
\lim_{N\to\infty}\frac{1}{N}\sum_{n=0}^{N-1}\int 
f\cdot T^n\E(f\mid\CZ_1)\cdot T^{2n}\E(f\mid\CZ_1)\, d\mu
\end{equation}
is positive.
Since $\CZ_1$ is $T$-invariant, $T^n\E(f\mid\CZ_1)\cdot T^{2n}\E(f\mid\CZ_1)$ 
is measurable with respect to $\CZ_1$ and so we can 
replace~\eqref{eq:intint} by
$$
\lim_{N\to\infty}\frac{1}{N}\sum_{n=0}^{N-1}\int 
\E(f\mid\CZ_1)\cdot T^n\E(f\mid\CZ_1)\cdot T^{2n}\E(f\mid\CZ_1)\, d\mu \ .
$$
This means that we can assume that the first term 
is also measurable with respect to the Kronecker factor, and so we can 
assume that 
$f$ is a nonnegative function that is measurable with 
respect to the Kronecker.  Thus the system $X$ can be assumed 
to be $Z_1$ and the transformation $T$ is rotation by some 
irrational $\alpha$.  Thus it suffices to show that 
$$\lim_{N\to\infty}\frac{1}{N}\sum_{n=0}^{N-1}\int_{Z_1}
f(s)\cdot f(s+n\alpha)\cdot f(s+2n\alpha)\,dm(s)$$ 
is positive.  
Since $\{n\alpha\}$ is equidistributed in $Z$, this limit 
approaches 
\begin{equation}\label{eq:iint}
\iint_{Z_1\times Z_1}f(s)\cdot f(s+t)\cdot f(s+2t)\, dm(s)dm(t) \ .
\end{equation}
But
$$
\lim_{t\to 0}\int_{Z_1}f(s)\cdot f(s+t)\cdot f(s+2t)\, dm(s) = 
\int_{Z_1}f(s)^3\, dm(s) \ ,
$$
which is clearly positive.  In particular, the double 
integral in~\eqref{eq:iint} is positive.  
\end{proof}

In the proof we have actually proven a stronger statement 
than needed to obtain Roth's Theorem: we have shown the 
existence of the limit of the double average in $L^2(\mu)$.  
Letting $\tilde{f}=\E(f\mid\CZ_1)$ for $f\in L^\infty(\mu)$, 
we have show that the double average~\eqref{eq:doubleav} 
converges to 
$$
\int_{Z_1}\tilde{f_1}(\pi_1(x) + s)\cdot 
\tilde{f_2}(\pi_1(x) + 2s)\, dm(s) \ .
$$

More generally, the same sort of argument can be used to show 
that in a weak mixing system, the Kronecker factor is {\em characteristic }
for the averages~\eqref{eq:APav} for all $k\geq 1$, 
meaning that to prove convergence of these average in a weak 
mixing system it suffices to assume that the system is a  
Kronecker system.  
Using Fourier analysis, one then gets convergence of the 
averages~\eqref{eq:APav} for weak mixing systems.

\subsection{Multiple averages}
\label{sec:higher}
We want to carry out similar analysis for the multiple averages
$$
\frac 1N\sum_{n=0}^{N-1} T^nf_1\cdot T^{2n}f_2\cdot\ldots\cdot 
T^{kn}f_k 
$$
and show the existence of the limit in $L^2(\mu)$ as $N\to\infty$.  
In his proof of Szemer\'edi's Theorem in~\cite{F} and subsequent 
proofs of Szemer\'edi's Theorem via ergodic theory such as~\cite{FKO}, 
the approach of Section~\ref{sec:k3} is not 
the one used for $k\geq 3$.  Namely, they do not show the 
existence of the limit and then analyze the limit itself to 
show it is positive.  A weaker statement is proved, only giving 
that the $\liminf$ of~\ref{eq:strongmultrec} is positive.  We 
will not discuss the intricate structure theorem 
and induction needed to prove this.

Already for convergence for 
$k=3$, one needs to consider more than just rotational behavior.  

\begin{example}
Given a system $(X, \CX, \mu, T)$, let $F(Tx) = f(x)F(x) \ ,
$
where
$$f(Tx) = \lambda f(x) \ \text{ and } \ |\lambda| = 1 \ .
$$
Then 
$$ F(T^{n}x) = f(x)f(Tx) \ldots f(T^{n-1}x)F(x)  = 
\lambda^{\frac{n(n-1)}{2}}\bigl(f(x)\bigr)^{n}F(x)
$$
and so 
$$
F(x) = \bigl(F(T^{n}x)\bigr)^{3}\bigl(F(T^{2n}x)\bigr)^{-3}F(T^{3n}x) \ .
$$
This means that there is some relation among
$$(x, T^{n}x, T^{2n}x, T^{3n}x)$$
that not arising from the Kronecker factor.  
\end{example}
One can construct more complicated examples (see Furstenberg~\cite{F6}) 
that show that even such generalized eigenfunctions do not suffice 
for determining the limiting behavior for $k=3$.  More precisely, 
the factor corresponding to generalized eigenfunctions (the {\em Abramov 
factor}) is not characteristic for the average~\ref{eq:APav} with $k=3$.

To understand the triple average, one needs to take into account 
systems more complicated than such Kronecker and Abramov 
systems.  The simplest such example
is a $2$-step nilsystem (the use of this terminology will be 
clarified later):
\begin{example}
Let $X = \T\times\T$, with Borel $\sigma$-algebra, and Haar measure.  
Fix $\alpha\in\T$ and define $T\colon X\to X$ by 
$$
T(x, y) = (x+\alpha, y+x) 
$$
The system is ergodic if and only if $\alpha\notin\Q$.  
\end{example}

The system is not isomorphic to a group rotation, as 
can be seen by defining $f(x,y) = e(y) = \exp(2\pi i y)$.  
Then for all $n\in\Z$, 
$$
T^n(x,y) = (x+n\alpha, y+nx+\frac{n(n-1)}{2}\alpha) 
$$
and so 
$$
f(T^n(x,y)) = e(y) e(nx)e\bigl(\frac{n(n-1)}{2}\alpha\bigr) \ .
$$
Quadratic expressions like these do not arise from a rotation 
on a group.

\section{The structure theorem}

\subsection{Major steps in the proof of Theorem~\ref{th:mainAP}}

In broad terms, there are four major steps in the proof of 
Theorem~\ref{th:mainAP}.  

For each $k\in\N$, we inductively define 
a seminorm $\nnorm{\cdot}_k$ that controls the asymptotic 
behavior of the average.  More precisely, we show that 
if $|f_1|\leq 1, \ldots, |f_k|\leq 1$, 
then 
\begin{equation}
\label{eq:charr}
\limsup_{N\to\infty}\Bigl\Vert\frac{1}{N}
\sum_{n=o}^{N-1}T^nf_1\cdot T^{2n}f_2\cdot\ldots\cdot 
T^{kn}f_k\Bigr\Vert_{L^2(\mu)}\leq \min_{1\leq j\leq k}\nnorm{f_j}_k \ .
\end{equation}

Using these seminorms, we 
define factors $Z_k$ of $X$ such that for 
$f\in L^\infty(\mu)$, 
$$\E(f\,\vert\,\CZ_{k-1}) = 0 \text{ if and only if } \nnorm f_k = 0 \ .
$$
It follows from~\ref{eq:charr} that 
the factor $Z_{k-1}$ is characteristic for the average~\eqref{eq:APav}.

The bulk of the work is then to give a 
``geometric'' description of these factors. 
This description is in terms of nilpotent groups, and more 
precisely we show that the dynamics of translation on 
homogeneous spaces of a nilpotent Lie group determines the 
limiting behavior of these averages.  
This is the content of the structure theorem, explained in 
Section~\ref{sec:structure}.  
(A more detailed expository version of this is given in 
Host~\cite{Host}; for full details, see~\cite{HK4}.)

Finally, we show convergence for these particular types of systems.  

Roughly speaking, this same outline applies to other convergence 
results we consider in the sequel, such as averages along polynomial 
times, averages 
along cubes, or averages for commuting transformations.  
For each average, we find a characteristic factor 
that can be described in geometric terms, allowing us to 
prove convergence in the characteristic factor.

\subsection{The role of nilsystems}
\label{sec:structure}
We have already seen that the limit behavior of the double average 
is controlled by group rotations, meaning the Kronecker 
factor is characteristic for this average.
Furthermore, we have seen that something more is 
needed to control the limit behavior of the triple average.  
Our goal here is to explain how the multiple averages of~\eqref{eq:APav}, 
and some more general averages, 
are controlled by nilsystems.  We start with some terminology.  

Let $G$ be a group. 
If $g,h\in G$, let $[g,h]=g\inv h\inv gh$ denote the commutator 
of $g$ and $h$.
If $A, B\subset G$, we write $[A,B]$ for the subgroup
of $G$ spanned by $\{[a,b]:a\in A,\ b\in B\}$. 
The lower central series
$$
G=G_1\supset G_2\supset \dots\supset G_j\supset G_{j+1}\supset\dots
$$
of $G$ is defined inductively, setting $G_1=G$ and 
$G_{j+1}=[G,G_j]$ for $j\geq 1$.  We say that 
$G$ is {\em $k$-step nilpotent} if $G_{k+1}=\{1\}$.

If $G$ is a $k$-step nilpotent Lie group and $\Gamma$ is a discrete
cocompact subgroup, the compact manifold $X=G/\Gamma$ is 
{\em $k$-step nilmanifold}.  

The group $G$ acts naturally on $X$ by
left translation:  if $a\in G$ and $x\in X$, the translation $T_a$ 
by $a$ is given by $T_a(x\Gamma) = (ax)\Gamma$.  
There is a unique Borel probability measure $\mu$ 
(the {\em Haar measure}) on $X$ that is invariant
under this action 
Fixing an element $a\in G$, 
the system $(G/\Gamma, \CG/\Gamma, T_a, \mu)$ 
is a {\em $k$-step nilsystem} and 
$T_a$ is a {\em nilrotation}.  

The system $(X,\CX, \mu, T)$ is an {\em inverse limit} of a
sequence of factors $\{(X_j,\CX_j, \mu_j, T)\}$ if $\{\CX_j\}_{j\in\N}$ is an
increasing sequence of $T$-invariant sub-$\sigma$-algebras 
such that $\bigvee_{j\in\N}\CX_j = \CX$ up to null
sets.
If each system $(X_{j}, \CX_j, \mu_j, T)$ 
is isomorphic to a $k$-step nilsystem, then 
$(X,\CX, \mu, T)$ is an {\em inverse limit of $k$-step nilsystems}.

Proving convergence of the averages~\eqref{eq:APav} is only possible 
if one can has a good description of some characteristic factor
for these averages.  This is the content of the structure theorem:
\begin{theorem}[Host and Kra~\cite{HK5}]
There exists a characteristic factor for the averages~\eqref{eq:APav} 
which is isomorphic to an inverse limit of $(k-1)$-step nilsystems.  
\end{theorem}

\subsection{Examples of nilsystems}
We give two examples of nilsystems that illustrate their general 
properties.
\begin{example}
Let $G=\Z\times\T\times \T$  with  multiplication given by 
$$
(k,x,y)*(k',x',y')=(k+k', x+x'\pmod 1, y+y'+ 2kx'\pmod 1)\ .
$$
The commutator subgroup of $G$ is $\{0\}\times \{0\}\times\T$, and so
$G$ is $2$-step nilpotent.
The subgroup $\Gamma=\Z\times\{0\}\times\{0\}$ is discrete and
cocompact, and thus $X = G/\Gamma$ is a nilmanifold.  
Let $\CX$ denote the Borel $\sigma$-algebra and let 
$\mu$ denote Haar measure on $X$.
Fix some irrational $\alpha\in\T$, let  $a=(1,\alpha,\alpha)$, 
and let $T:X\to X$ be translation by $a$.
Then $(X,\mu,T)$ is a $2$-step nilsystem.

The Kronecker factor of $X$ is $\T$ with rotation by $\alpha$.
Identifying $X$ with $\T^2$ via the map 
$(k,x,y)\mapsto (x,y)$, the transformation $T$  takes on 
the familiar form of a skew transformation:
$$
T(x,y) = (x+\alpha, y+2x+\alpha) \ .
$$
This system is ergodic if and only if $\alpha\notin\Q$:  for 
$x,y\in X$ and $n\in\Z$, 
$$
T^n(x,y) = (x+n\alpha, y+2nx+n^2\alpha)
$$
and equidistribution of the sequence $\{T^n(x,y)\}$ is equivalent 
to ergodicity.  
\end{example}

\begin{example}
Let $G$ be the Heisenberg group $\R\times\R\times\R$ with
multiplication given by
$$
(x,y,z)*(x',y',z') = (x+x', y+y', z+z'+xy') \ .
$$
Then $G$ is a $2$-step nilpotent Lie group.
The subgroup $\Gamma = \Z\times\Z\times\Z$ is discrete and cocompact and 
so  $X=G/\Gamma$ is a nilmanifold.  Letting
$T$ be the translation by $a = (a_1, a_2, a_3)\in G$ where 
$a_1, a_2$ are independent over $\Q$ and $a_3\in\R$, and taking $\CX$ 
to be the Borel $\sigma$-algebra and $\mu$ the Haar measure, we have 
that $(X,\CX, \mu, T)$ is a nilsystem.  The system is 
ergodic if and only if $a_1, a_2$ are independent over 
$\Q$

The compact abelian group $G/G_2\Gamma$ is 
isomorphic to $\T^{2}$ and the rotation on $\T^2$ 
by $(a_1, a_2)$ is ergodic.  
The Kronecker factor 
of $X$ is the factor induced by functions on $x_1, x_2$.  

The system $(X, \CX, \mu, T)$ is (uniquely) ergodic.
\end{example}

The dynamics of the first example gives rise to 
quadratic sequences, such as $\{n^2\alpha\}$, and the 
dynamics of the second example gives rise to generalized 
quadratic sequences such as $\{\lfloor n\alpha\rfloor n\beta\}$.

\subsection{Motivation for nilpotent groups}
The content of the Structure Theorem is that 
nilpotent groups, or more precisely the dynamics 
of a translation on the homogeneous space of a nilpotent Lie group, 
control the limiting behavior of the averages along arithmetic 
progressions.  We give some motivation as to why nilpotent 
groups arise.  

If $G$ is an abelian group, then 
$$
\{(g, gz, gz^{2}, \ldots, gz^n)\colon g, z\in G\}
$$ 
is a subgroup of $G^n$.  
However, this does not hold if $G$ is not abelian.  
To make these arithmetic progressions into a group, one 
must take into account the commutators.  This 
is the content of the following theorem, proven in 
different contexts by Hall~\cite{H}, Petresco~\cite{P}, Lazard~\cite{La}, 
Leibman~\cite{L}: 

\begin{theorem}
If $G$ is a group, then for any $x,y\in G$, there exist $z\in G$ and 
$w_i\in G_i$ such that 
\begin{multline*}
(x, x^2, x^3, \ldots, x^n)\times (y, y^2, y^3, \ldots, y^n) = \\
(z, z^2w_1, z^3w_1^3w_2,\ldots, z_{\ }^{\binom{n}{1}}
w_1^{\binom{n}{2}}
w_2^{\binom{n}{3}}\ldots 
w_{n-1}^{\binom{n}{n}}) \ .
\end{multline*}
Furthermore, these expressions form a group.  
\end{theorem}

If $G$ is a group, a {\em geometric progression} is a sequence of the form 
$$
(g, gz, gz^2w_1, gz^3w_1^3w_2, \ldots, gz^{\binom{n}{1}}
w_1^{\binom{n}{2}}
\ldots 
w_{n-1}^{\binom{n}{n}}) \ ,
$$
where $g,z\in G$ and $w_i\in G_i$.   

Thus if $G$ is abelian, $g$ and $z$ determine the whole sequence.  
On the other hand, if $G$ is $k$-step nilpotent with $k<n$, 
the first $k$ terms determine the whole sequence.  

Similarly, if $(G/\Gamma, \CG/\Gamma, \mu, T_a)$ is a $k$-step nilsystem and 
$$
x_1 = g_1\Gamma, x_2=g_2\Gamma, \ldots, x_k=g_k\Gamma, \ldots, x_n=g_n\Gamma
$$ 
is a geometric progression, then the first $k$ terms determine 
the rest.  
Thus $a^{k+1}x\Gamma$ is a function of the first $k$ terms 
$ax\Gamma, a^2x\Gamma, \ldots, a^{k}x\Gamma $.  

This means that the $(k+1)$-st term $T^{(k+1)n}x$ in an arithmetic 
progression $T^nx, \ldots,$ $T^{kn}x$ is
constrained by first $k$ terms.  
More interestingly, the converse also holds: 
in an arbitrary system $(X,\CX, \mu,T)$, any $k$-step nilpotent factor 
places a constraint on $(x, T^nx, T^{2n}x, \ldots, T^{kn}x)$.

\section{Building characteristic factors}
\label{sec:charfac}

The material in this and the next section is based on~\cite{HK4} 
and the reader is referred to~\cite{HK4} for full proofs.  
To describe characteristic factors for the averages~\eqref{eq:APav}, 
for each $k\in\N$ we define a seminorm and use it to define
these factors.  We start by defining certain 
measures that are then used to define the seminorms.
Throughout this section, we assume that $(X, \CX, \mu, T)$ is 
an ergodic system. 

\subsection{Definition of the measures}
Let $X\type k= X^{2^k}$ and define 
$T\type k\colon X\type k\to X\type k$ by 
$T\type k=T\times\dots\times T$ (taken $2^k$ 
times).

We write a point $\bx\in X\type k$ as 
$\bx=\bigl(x_\epsilon:\epsilon\in\{0,1\}^k\bigr)$ and 
make the natural identification of 
$X\type{k+1}$ with $X\type k\times X\type k$, 
writing $\bx=(\bx',\bx'')$
for a point of $X\type{k+1}$, with $\bx',\bx''\in X\type k$. 

By induction, we define a measure $\mu\type k$ on $X\type k$ 
invariant under $T\type k$.  
Set $\mu\type 0:=\mu$.
Let $\CI\type k$ be the invariant $\sigma$-algebra of
$(X\type k,\CX\type k, \mu\type k,T\type k)$.  (Note that this system 
is not necessarily ergodic.)
Then $\mu\type{k+1}$ is defined to be 
the {\em relatively independent joining} 
of $\mu\type k$ with itself over $\CI\type k$, meaning that
if $F$ and $G$ are  bounded functions on $X\type k$,
\begin{equation}
\label{eq:defmuk}
 \int_{X\type{k+1}}F(\bx')\cdot G(\bx'')\,d\mu\type{k+1}(\bx)
= 
\int_{X\type k}\E(F\mid\CI\type k)(\by)\cdot\E(G\mid\CI\type k)(\by)
\,d\mu\type k(\by) \ .
\end{equation}

Since $(X, \CX, \mu, T)$ is assumed to be ergodic, $\CI\type 0$ 
is trivial and $\mu\type 1 = \mu\times \mu$.  If the system 
is weak mixing, then for all $k\geq 1$, $\mu\type k$ is the product measure 
$\mu\times\mu\times\ldots\times\mu$, taken $2^k$ times.

\subsection{Symmetries of the measures}

Writing a point $\bx\in X\type k$ as 
$$
\bx = (x_\epsilon\colon \epsilon\in\{0,1\}^k) \ ,
$$
we identify the indexing set $\{0,1\}^k$ of this point 
with the vertices of the Euclidean cube.  

An isometry $\sigma$ of $\{0,1\}^k$ induces a map $\sigma_*\colon 
X\type k\to X\type k$ by {\em permuting the coordinates}:
$$(\sigma_*(\bx))_\epsilon = x_{\sigma(\epsilon)}\ .$$

For example, from the diagonal symmetries for $k=2$, we have 
the permutations
$$
(x_{00}, x_{01}, x_{10}, x_{11})\mapsto (x_{00}, x_{10}, x_{01}, x_{11})\ \ 
$$
$$
(x_{00}, x_{01}, x_{10}, x_{11})\mapsto (x_{11}, x_{01}, x_{10}, x_{00})\ .
$$

By induction, the measures are invariant under permutations:
\begin{lemma}
\label{lemma:symmetry}
For each $k\in\N$, 
the measure $\mu\type k$ is invariant under all permutations of 
coordinates arising from isometries of the unit Euclidean cube.
\end{lemma}

\subsection{Defining seminorms}

For each $k\in\N$, we define a 
seminorm on $L^\infty(\mu)$ by setting 
$$
\nnorm f_k^{2^k}= \int_{X\type k}\prod_{\epsilon\in\{0,1\}^k}f(x_\epsilon)
\,d\mu\type k(\bx) \ .
$$

By definition of the measure $\mu\type k$, this integral is 
equal to 
$$
\int_{X\type {k-1}}\E\Bigl(\prod_{\epsilon\in\{0,1\}^{k-1}}f(x_\epsilon)
\mid\CI\type{k-1}\Bigr)^2
\,d\mu\type {k-1}
$$
and so in particular it is nonnegative.

From the symmetries of the measure $\mu\type k$ (Lemma~\ref{lemma:symmetry}), 
we have a version of the Cauchy-Schwarz inequality for the seminorms, 
referred to as a Cauchy-Schwarz-Gowers inequality:
\begin{lemma}
For $\epsilon\in\{0,1\}^k$, let $f_\epsilon\in L^\infty(\mu)$.  Then 
$$
\Bigl\vert \int\prod_{\epsilon\in\{0,1\}^k}f_\epsilon(x_\epsilon)
\,d\mu\type k(\bx)\Bigr\vert \leq 
\prod_{\epsilon\in\{0,1\}^k}\nnorm f_k \ .
$$
\end{lemma}

As a corollary, the map $f\mapsto\nnorm f_k$ is subadditive and so:

\begin{corollary} 
For every $k\in\N$, 
$\nnorm\cdot_k$ is a seminorm on $L^\infty(\mu)$.
\end{corollary}

Since the system $(X, \CX, \mu, T)$ is ergodic,  the 
$\sigma$-algebra $\CI\type 0$ is trivial, $\mu\type 1 = \mu\times\mu$ 
and $\nnorm f_1=\Bigl|\int f\,d\mu\Bigr|$.
By induction,
$$\nnorm f_1\leq\nnorm f_2\leq\cdots\leq\nnorm f_k
\leq\cdots\leq\norm f_\infty\ .$$
If the system is weak mixing, then $\nnorm f_k = \nnorm f_1$ for 
all $k\in\N$.

By induction and the ergodic theorem, we have a second 
presentation of these seminorms:
\begin{lemma}
For every $k\geq 1$, 
$$
\nnorm f_{k+1}^{2^{k+1}}
=\lim_{N\to\infty} \frac 1N \sum_{n=0}^{N-1}
\nnorm{f\cdot T^nf}_k^{2^{k}}
\ .
$$
\end{lemma}

\subsection{Seminorms control the averages~\ref{eq:APav}}
The seminorms $\nnorm\cdot_k$ 
control the averages along arithmetic progressions:
\begin{lemma}
\label{lemma:charseminorm}
Assume that $(X, \CX, \mu, T)$ is ergodic and let $k\in\N$.  
If $\norm{f_1}_\infty$,\dots, $\norm{f_k}_\infty \leq 1$, 
then 
$$
\limsup_{N\to\infty}\Bigl\Vert \frac 1N \sum_{n=0}^{N-1}
 T^nf_1\cdot T^{2n}f_2\cdot\ldots \cdot T^{kn}f_k
 \Bigr\Vert_2  \leq 
\min_{\ell=1, \ldots, k}\ell\nnorm{f_\ell}_k\ .
$$
\end{lemma}

\begin{proof}
We proceed by induction on $k$.  For $k=1$ this is 
trivial.  Assume it holds for $k\geq 1$.
Define $u_n=T^nf_1\cdot T^{2n}f_2\cdots 
T^{(k+1)n} f_{k+1}$, and assume that $\ell>1$ (the case $\ell=1$ is 
similar). Then
\begin{multline*}
\bigl|\frac 1N\sum_{n=0}^{N-1} \langle u_{n+h},u_n\rangle\Bigr|=
\Bigl| \int (f_1\cdot T^hf_1)\frac 1N\sum_{n=0}^{N-1}
\prod_{j=2}^{k+1} T^{(j-1)n}(f_j\cdot T^{jh}f_j)\,d\mu\Bigr|\\
\leq \Bigl\Vert \frac 1N\sum_{n=0}^{N-1}
\prod_{j=2}^{k+1} T^{(j-1)n}(f_j\cdot T^{jh}f_j)\Bigr\Vert_2\ .
\end{multline*}

By the induction hypothesis, 
$\gamma_h  \leq \ell \nnorm{f_\ell\cdot T^{\ell h}f_\ell}_k$.
Thus 
 $$
 \frac 1H\sum_{h=0}^{H-1}\gamma_h\leq \ell\,\frac \ell{\ell H}
 \sum_{n=0}^{\ell H-1}\nnorm{f_\ell\cdot T^nf_\ell}_k
$$
and the statement follows from the van der Corput Lemma 
(Lemma~\ref{lem:vdc}) and the definition of the seminorm
 $\nnorm\cdot_{k+1}$.
\end{proof}

\subsection{The Kronecker factor, revisited ($k=2$)}

We have seen two presentations of the Kronecker factor $(Z_1, \CZ_1, m, T)$: 
it is the largest abelian group rotation factor and it is 
the sub-$\sigma$-algebra of $\CX$ generated by the 
eigenfunctions.   Another equivalent formulation is that it is the smallest 
sub-$\sigma$-algebra of $\CX$ such that all invariant functions 
of $(X\times X, \CX\times\CX, \mu\times\mu, T\times T)$ are 
measurable with respect to $\CZ_1\times\CZ_1$.  Recall that 
$\pi_1\colon X\to\CZ_1$ denotes the factor map.  

We give an explicit description of the measure $\mu\type 2$, 
and thus give yet another description of the Kronecker factor.  
For $f\in L^\infty(\mu)$, write $\tilde{f}=\E(f\mid \CZ_1)$.

For $s\in Z_1$ and $f_0, f_1\in L^\infty(\mu)$, 
we define a probability measure $\mu_s$ on $X\times X$ by 
$$\int_{X\times X}f_0(x_0) f_1(x_1)\,d\mu_s(x_0,x_1):=
\int_{Z_1}\tf_0(z)\tilde f_1(z+s)\,dm(z) \ .
$$
This measure is $T\times T$-invariant and the ergodic 
decomposition of $\mu\times\mu$ under $T\times T$ is given by 
$$
\mu\times\mu=\int_{Z_1}\mu_s\,dm(s)\ .
$$
Thus for $m$-almost every $s\in Z_1$, the system $(X\times X, \CX\times\CX, 
\mu_s, T\times T)$ is ergodic and 
$$
\mu\type 2=\int_{Z_1}\mu_s\times\mu_s\,dm(s) \ .
$$
More generally, if $f_\epsilon$, $\epsilon\in\{0,1\}^2$, are 
measurable functions on $X$, then 
\begin{multline*}
\int_{X\type 2} f_{00}\otimes f_{01}\otimes f_{10}\otimes f_{11}
\,d\mu\type 2 \\
=\int_{Z_1^3}\tf_{00}(z)\cdot\tf_{01}(z+s)\cdot\tf_{10}(z+t)\cdot
\tf_{11}(z+s+t)\,dm(z)\,dm(s)\,dm(t) \ .
\end{multline*}
It follows immediately that:
\begin{multline*} 
\nnorm f_2^4
:=\int f\otimes f\otimes f\otimes f\,d\mu\type 2 \\
=\int_{Z_1^3}\tf(z)\cdot\tf(z+s)\cdot
\tf(z+t)\cdot\tf(z+s+t)\,dm(z)\,dm(s)\,dm(t) \ .
\end{multline*}

As a corollary, $\nnorm f_2$ is the $\ell^4$-norm of the Fourier Transform 
of $\tf$ and the factor $Z_1$, 
defined by $\nnorm f_2 = 0$ if and only if $\E(f\mid\CZ_1) = 0$ 
for $f\in L^\infty(\mu)$, 
 is the Kronecker factor of $(X, \CX, \mu, T)$.

\subsection{Factors for all $k\geq 1$}

Using these seminorms, 
we define factors $Z_k=Z_k(X)$ for $k\geq 1$ of $X$ that 
generalize the relation between the Kronecker factor 
$Z_1$ and the second seminorm $\nnorm \cdot_2$: 
for $f\in L^\infty(\mu)$, $\E(f\mid \CZ_k)=0$ if
and only if $\nnorm f_{k+1}=0$.  
To explain this, we start by describing some  geometric 
properties of the measures $\mu\type k$.

Indexing $X\type k$ by the coordinates $\{0,1\}^k$ 
of the Euclidean cube, it is natural to use geometric terms like 
{\em side, edge, vertex} for subsets of $\{0,1\}^k$.  
For example, the following illustrates 
the point $\bx\in X\type 3$ with 
the side $\alpha=\{010, 011, 110, 111\}$:

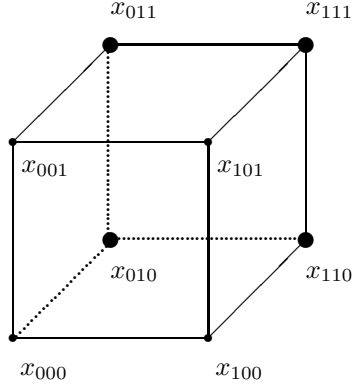
\begin{center}
\setlength{\unitlength}{1.3cm}
\begin{picture}(3.5,3.5)(-0.5,0)
\put(1,3){\circle*{0.15}} 
\put(3,3){\circle*{0.15}} 
\put(3,1){\circle*{0.15}} 
\put(2,2){\circle*{0.08}} 
\put(0,2){\circle*{0.08}} 
\put(0,0){\circle*{0.08}} 
\put(2,0){\circle*{0.08}} 
\put(1,1){\circle*{0.15}} 
\thinlines
\put(2,0) {\line(0,1){2}}
\put(2,0) {\line(1,1){1}}
\put(2,2) {\line(1,1){1}}
\put(3,1) {\line(0,1){2}}
\put(0,0) {\line(1,0){2}}
\put(0,0) {\line(0,1){2}}
\put(1,3) {\line(1,0){2}}
\put(0,2) {\line(1,0){2}}
\put(0,2){\line(1,1){1}}
\multiput(0,0)(0.04,0.04){25}{.}
\multiput(1,1)(0.05,0){40}{.}
\multiput(1,1)(0,0.05){40}{\kern -0.8mm .}

\put(1,3){\raisebox{4mm}{$x_{011}$}}
\put(3,3){\raisebox{4mm}{$x_{111}$}}
\put(3,1){\raisebox{-5mm}{$x_{110}$}}
\put(2,2){\raisebox{-4mm}{\ $x_{101}$}}
\put(0,2){\raisebox{-4mm}{\ $x_{001}$}}
\put(0,0){\raisebox{-5mm}{$\;x_{000}$}}
\put(2,0){\raisebox{-5mm}{$\;x_{100}$}}
\put(1,1){\raisebox{-5mm}{$x_{010}$}}
\end{picture}\\
\end{center}

\bigskip
\medskip

Let $\alpha\subset\{0,1\}^k$ be a side.
The {\em side transformation} $T\type k_\alpha$
of  $X\type k$ is defined by:
$$
 \bigl(T\type k_\alpha\bx\bigr)_\beps=
\begin{cases}
Tx_\beps & \text{ if }\beps\in\alpha\ ;\\
x_\beps&\text{ otherwise \ .}
\end{cases}
$$

We can represent the transformation $T_\alpha$ associated 
to the side $\{010, 011, 110, 111\}$ by:
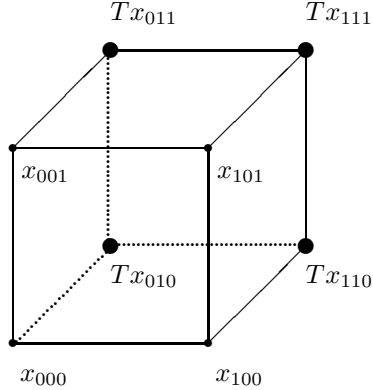
\begin{center}
\setlength{\unitlength}{1.3cm}
\begin{picture}(3.5,3.5)(-0.5,0)
\put(1,3){\circle*{0.15}} 
\put(3,3){\circle*{0.15}} 
\put(3,1){\circle*{0.15}} 
\put(2,2){\circle*{0.08}} 
\put(0,2){\circle*{0.08}} 
\put(0,0){\circle*{0.08}} 
\put(2,0){\circle*{0.08}} 
\put(1,1){\circle*{0.15}} 
\thinlines
\put(2,0) {\line(0,1){2}}
\put(2,0) {\line(1,1){1}}
\put(2,2) {\line(1,1){1}}
\put(3,1) {\line(0,1){2}}
\put(0,0) {\line(1,0){2}}
\put(0,0) {\line(0,1){2}}
\put(1,3) {\line(1,0){2}}
\put(0,2) {\line(1,0){2}}
\put(0,2){\line(1,1){1}}
\multiput(0,0)(0.04,0.04){25}{.}
\multiput(1,1)(0.05,0){40}{.}
\multiput(1,1)(0,0.05){40}{\kern -0.8mm .}

\put(1,3){\raisebox{4mm}{$Tx_{011}$}}
\put(3,3){\raisebox{4mm}{$Tx_{111}$}}
\put(3,1){\raisebox{-5mm}{$Tx_{110}$}}
\put(2,2){\raisebox{-4mm}{\ $x_{101}$}}
\put(0,2){\raisebox{-4mm}{\ $x_{001}$}}
\put(0,0){\raisebox{-5mm}{$\;x_{000}$}}
\put(2,0){\raisebox{-5mm}{$\;x_{100}$}}
\put(1,1){\raisebox{-5mm}{$Tx_{010}$}}
\end{picture}
\end{center}

\bigskip
\medskip

Since permutations of coordinates leave the measure $\mu\type k$ 
invariant and act transitively on the sides, we have:
\begin{lemma}
For all $k\in\N$, the 
measure $\mu\type k$ is invariant under the side transformations.
\end{lemma}

We now view $X\type k$ in a different way, identifying 
$X\type k = X\times X^{2^k-1}$.    A point $\bx\in X\type k$ 
is now written as 
$$
\bx = (x_\bzero, \tilde{x}) \text{ where }\tilde{x} \in 
X^{2^k-1}, x_\bzero\in X, \text{ and } \bzero = (00\ldots 0)\in\{0,1\}^k \ .
$$
Although the $0$ coordinate has been singled out and seems 
to play a particular role, it follows from the symmetries of the 
measure $\mu\type k$ (Lemma~\ref{lemma:symmetry}) that any 
other coordinate could have been used instead.  

If $\alpha\subset \{0,1\}^k$ is a side that does not contain $\bzero$ 
(there are $k$ such sides), the transformation 
$T_\alpha\type k$ leaves the coordinate $\bzero$ invariant.  
It follows from induction and the definition of the measure $\mu\type k$ 
that:

\begin{proposition} 
Let $k\in\N$.  If $B\subset X^{2^k-1}$, there exists $A\subset X$  with 
\begin{equation}
\label{eq:deffactor}
\one_A(x_\bzero) = \one_B(\tilde{x}) \text{ for almost 
all } \bx\in X\type k\end{equation}
if and only if $X\times B$ is invariant under the $k$ transformations 
$T_\alpha\type k$ arising from the $k$ sides $\alpha$ 
not containing $\bzero$.  
\end{proposition}

This means that 
the subsets $A\subset X$ such that there exists $B\subset X^{2^k-1}$ 
satisfying~\eqref{eq:deffactor}
form an invariant sub-$\sigma$-algebra $\CZ_{k-1} = \CZ_{k-1}(X)$ 
of $\CX$.  We define $Z_{k-1} = Z_{k-1}(X)$ to be the associated factor.  
Thus $\CZ_{k-1}(X)$ is defined to be the sub-$\sigma$-algebra of 
sets $A\subset X$ such that Equation~\eqref{eq:deffactor} 
holds for some set $B\subset X^{2^k-1}$.  

We give some properties of the factors:
\begin{proposition}
\begin{enumerate}
\item 
For every bounded function $f$ on $X$, 
$$
\nnorm f_k = 0 \text{ if and only if } \E(f\mid\CZ_{k-1}) = 0 \ .
$$

\item
For bounded functions $f_\epsilon$, $\epsilon\in\{0,1\}^k$, on $X$, 
$$
\int\prod_{\epsilon\in\{0,1\}^k}f_\epsilon(x_\epsilon)\, d\mu\type k(\bx) 
 = \int\prod_{\epsilon\in\{0,1\}^k}
\E(f_\epsilon\mid\CZ_{k-1})(x_\epsilon)\, d\mu\type k(\bx) \ .
$$
Furthermore, $\CZ_{k-1}$ is the smallest sub-$\sigma$-algebra of $\CX$ with 
this property.  

\item The invariant sets of $(X\type k, \CX\type k, \mu\type k, T\type k)$ 
are measurable with respect to $\CZ_k\type k$.  Furthermore, 
$\CZ_k$ is the smallest sub-$\sigma$-algebra of $\CX$ with this property.  
\end{enumerate}
\end{proposition}

The proof of this proposition relies on showing a similar 
formula to that used (in Equation~\eqref{eq:defmuk}) to define 
the measures $\mu\type k$, but with respect to the new 
identification separating $1$ coordinate from the $2^k-1$ others.  
Namely, for 
bounded functions $f$ on $X$ and $F$ on $X^{2^k-1}$, 
$$
\int_{X\type k} f(x_{\bf 0})\cdot F(\tilde{x})\, d\mu\type k(\bx) = 
\int_{X\type{k-1}}\E(f\mid\CZ_{k-1})\cdot\E(F\mid \CZ_{k-1})\, d\mu\type{k-1} \ .
$$
The given properties then follow using induction and the symmetries 
of the measures.

We have already seen that 
$Z_0$ is the trivial factor and $Z_1$ is the Kronecker factor.  
More generally, the sequence of factors is increasing:
$$
 Z_0\leftarrow Z_1\leftarrow\dots\leftarrow Z_k\leftarrow 
 Z_{k+1}\leftarrow\dots\leftarrow X\ .
$$
If $X$ is weak mixing, then $Z_k(X)$ is the trivial factor for every $k$.

An immediate consequence of Lemma~\ref{lemma:charseminorm} 
and the definition of the factors is that 
the factor $Z_{k-1}$ is characteristic for the average along 
arithmetic progressions:
\begin{proposition}
\label{prop:charfac}
For all $k\geq 1$, the factor $Z_{k-1}$ is 
characteristic for the convergence of the 
averages 
$$
\frac 1N \sum_{n=0}^{N-1}
 T^nf_1\cdot T^{2n}f_2\cdot\ldots\cdot T^{kn}f_k \ .
$$
\end{proposition}

\section{Structure theorem}
\subsection{Systems of order k}

For $k\geq 0$, 
an ergodic system $X$ is said to be of {\em order $k$} if $Z_k(X) = X$.  
This means that 
$\nnorm{\cdot}_{k+1}$ is a norm on $L^\infty(\mu)$.

Given an ergodic system $(X, \CX, \mu, T)$, 
$Z_k(X)$ is a system of order $k$, since $Z_k(Z_k(X)) = Z_k(X)$.  
The unique system of order zero is 
the trivial system, and a system of order $1$ is an ergodic rotation.
By definition, if a system is of order $k$, then it is also 
of order $k'$ for any $k'> k$.

By Proposition~\ref{prop:charfac}, to show convergence of 
$$
\frac{1}{N}\sum_{n=0}^{N-1}T^nf_1\cdot T^{2n}f_2\cdot\ldots\cdot T^{kn}f_k
$$
on an arbitrary system, it suffices to assume that each function 
is defined on the factor $Z_{k-1}$.  But 
since $Z_{k-1}(X)$ is a system of order $k$, it 
suffices to prove convergence of this average for systems 
of order $k-1$.  

In this language, the structure theorem becomes:
\begin{theorem}[Host and Kra~\cite{HK5}] 
A system of order $k$ is the inverse limit of a sequence of $k$-step 
nilsystems.
\end{theorem}

Before turning to the proof of the structure theorem, 
we show convergence for the average along arithmetic progressions 
in a nilsystem.

\subsection{Convergence on a nilmanifold}
Using general properties of nilmanifolds (see Furstenberg~\cite{F3} and 
Parry~\cite{Parry}), Lesigne~\cite{Lesigne} showed 
for connected group $G$ and Leibman~\cite{Leibman2} showed 
in the general case, convergence 
in a nilsystem:
\begin{theorem}
\label{th:convergence}
If $(X = G/\Gamma, \CG/\Gamma, \mu, T)$ is a nilsystem and $f$ 
is a continuous function on $X$, then 
$$
\frac{1}{N}\sum_{n=0}^{N-1} f(T^nx)
$$ 
converges for every $x\in X$.
\end{theorem}
(See also Ratner~\cite{Ratner} and Shah~\cite{Shah} for 
related convergence results.)

As a corollary, we have convergence in $L^2(\mu)$ for the 
average along arithmetic progressions in a nilmanifold:
\begin{corollary}
If $(X = G/\Gamma, \CG/\Gamma, \mu, T)$ is a nilsystem, 
$k\in\N$, and $f_1, f_2, \ldots, f_k\in L^\infty(\mu)$, then  
$$\lim_{N\to\infty}\frac{1}{N}\sum_{n=0}^{N-1}
T^nf_1\cdot T^{2n}f_2\cdot\ldots\cdot 
T^{kn}f_k
$$ 
exists in $L^2(\mu)$.  
\end{corollary}

\begin{proof}  By density, we can assume that the functions are continuous.  
By assumption, 
$G^k$ is a nilpotent Lie group, 
$\Gamma^k$ is a discrete cocompact subgroup and $X^k = G^k/\Gamma^k$ 
is a nilmanifold.  Let 
$$
s = (t, t^2, \ldots, t^k)\in G^k
$$
and let $S\colon X^k\to X^k$ be the translation by $s$, meaning that 
$$
S  = T\times T^2\times\ldots\times T^k \ .
$$
We apply Theorem~\ref{th:convergence} to 
$(X^k, S)$ with the continuous function 
$$
F(x_1, x_2, \ldots, x_k) = f_1(x_1)f_2(x_2)\ldots f_k(x_k)$$ 
at the point $y = (x, x, \ldots, x)$ and so the averages 
converge everywhere.
\end{proof}

Thus Theorem~\ref{th:mainAP} holds in a nilsystem, and we 
are left with proving the Structure Theorem.  

\subsection{A group of transformations}
To each ergodic system, we associate a group of 
measure preserving transformations.  The general 
approach is to show that for sufficiently many systems of 
order $k$, this group is a nilpotent Lie group.  The bulk of 
the work is to then show that this group 
acts transitively on the system.  Thus the system 
can be given the 
structure of a nilmanifold and the Structure Theorem follows.  

Most proofs are sketched or omitted completely, 
and the reader is referred to~\cite{HK5} 
for the details.

Let $(X, \CX, \mu, T)$ be an ergodic system.  
If $S\colon X\to X$  and
$\alpha\subset\{0,1\}^k$, define $S\type k_\alpha\colon X\type k 
\to X\type k$ by:
$$
 \bigl(S\type k_\alpha\bx\bigr)_\beps=
\begin{cases}
Sx_\beps & \text{if }\beps\in\alpha;\\
x_\beps & \text{otherwise\ .}
\end{cases}
$$

Let $\CG=\CG(X)$ be 
the group of transformations $S\colon X\to X$ such
that for all $k\in\N$ and all sides $\alpha\subset\{0,1\}^k$, the
measure $\mu\type k$ is invariant under $S\type k_\alpha$.

Some properties of this group are immediate.  
By symmetry, it suffices to consider one side.  
By definition, $T\in\CG$, and if $ST=TS$ then we also 
have that $S\in\CG$.  
If $S\in\CG$ and $k\in N$, then 
$\mu\type k$ is invariant under 
$S\type k\colon X\type k\to X\type k$.
Furthermore, $S\type kE=E$ for every $E\in\CI\type k$.

By induction, the invariance of the measure $\mu\type k$ under 
the side transformations, and commutator relations, we 
have:
\begin{proposition}
If $X$ is a system of order $k$, then $\CG(X)$ is a $k$-step 
nilpotent group.  
\end{proposition}

\subsection{Proof of the structure theorem}

We proceed by induction.  By the inductive assumption, we 
can assume that we are given a system $(X,\CX, \mu,T)$ of order $k$.  
We have a factor $(Y,\CY,\nu,T)$, where 
$Y=Z_{k-1}(X)$ and $\pi\colon X\to Y$ is the factor map.
Furthermore, $Y$ is an inverse limit of a sequence of $(k-1)$-step
nilsystems
$$
 Y=\limproj Y_i\ ;\ Y_i=G_i/\Gamma_i \ .
$$
We want to show that $X$ is an inverse limit of $k$-step nilsystems.

We have already shown that
if $f_\beps$, $\beps\in\{0,1\}^k$, are bounded functions on $X$, 
then 
$$
\int\prod_{\beps\in\{0,1\}^k}f_\beps(x_\beps)\,d\mu\type k(\bx)
=\int\prod_{\beps\in\{0,1\}^k}
\E\bigl(f_\beps\mid \CY)(x_\beps)\,d\mu\type k(\bx)
$$
In particular, for $f\in L^\infty(\mu)$,
$$
 \nnorm f_k=0 \text{ if and only if } \E(f\mid \CY)=0\ .
$$
Furthermore, 
$X$ does not admit a strict sub-$\sigma$-algebra $\CZ$ such that all 
invariant sets of $(X\type k,\mu\type k,T\type k)$ are measurable
with respect to $\CZ\type k$.
Recall also that the system  $(X\type k,\mu\type k,T\type k)$ is 
defined as a relatively independent
joining.

In ~\cite{F}, Furstenberg described the invariant $\sigma$-algebra 
for relatively independent joinings.  
It follows that $X$ is an {\em isometric extension} of $Y$, 
meaning that 
$X=Y\times H/K$ where $H$ is a compact group and $K$ is a closed subgroup, 
$\mu=\nu\times m$, where $m$ is the Haar measure of $H/K$, and 
the transformation $T$ is given by
$$
 T(y,u)=(Ty,\rho(y)\cdot u)
$$
for some map $\rho\colon Y\to H$.

\begin{lemma}
For every $h\in H$, the transformation  $(y,u)\mapsto (y,h\cdot u)$
of $X$ belongs to the center of $\CG(X)$.
\end{lemma}

Thus $H$ is abelian.  We can substitute $H/K$ for $H$, and 
we use additive notation for $H$.

We therefore have more information: $X$ is an {\em abelian extension} of $Y$, 
meaning that 
$X=Y\times H$ for some compact abelian group $H$, 
$\mu=\nu\times m$, where $m$ is the Haar measure of $H$, 
and the transformation $T$ is given by 
$T(y,u)=(Ty,u+\rho(y))$
for some map $\rho\colon Y\to H$.  We 
call $\rho$ the {\em cocycle} defining the
extension.

Furthermore, we show that the cocycle defining this 
extension has a particular form:
\begin{proposition}[The functional equation]
\label{prop:funceq}
If $(X, \CX, \mu, T)$ is a system of order $k$ and 
$(Y, \CY, \nu, T) = Z_{k-1}(X)$, then $X$ is an 
abelian extension of $Y$ via a compact group $H$ and 
for the cocycle $\rho$ defining this extension, 
there exists a map $\Phi\colon Y\type k\to H$ such that
\begin{equation}
\label{eq:funceq}
 \sum_{\beps\in\{0,1\}^k}(-1)^{\epsilon_1+\dots+\epsilon_k}
\rho(y_\beps)=\Phi(T\type k\by)-\Phi(\by)
\end{equation}
for $\nu\type k$-a.e. $\by\in\ Y\type k$.  
\end{proposition}

We can make a few more assumptions on our system. Namely, 
by induction we can deduce that 
$H$ is connected. 
Since every connected compact abelian group $H$ is an inverse limit of a
sequence of tori, we can further reduce to the case that $H=\T^d$.

\subsection{The case $k=2$ (The Conze-Lesigne Equation)}

We maintain notation of the preceding 
section and review what this means for the case $k=2$.  
By assumption, we have that 
$(Y,\CY, \nu,T)$ is a system of order $1$, meaning it is a group rotation.
The measure $\nu\type 2$ is the Haar measure of the subgroup
$$
\bigl\{ (y,y+s,y+t,y+s+t)\colon y,s,t\in Y\bigr\}
$$
of $Y^4$.
The functional equation of 
Proposition~\ref{prop:funceq} is: 
there exists $\Phi\colon Y^3\to \T^d$ with
$$ 
\rho(y)-\rho(y+s)-\rho(y+t)+\rho(y+s+t)
=\Phi(y+1,s,t)-\Phi(y,s,t)
$$

It follows that for every $s\in Y$, there exists $\phi_s\colon Y\to 
\T^d$ 
and $c_s\in \T^d$ satisfying the {\em Conze-Lesigne Equation} 
(see~\cite{CL4}):
\begin{equation}
\label{eq:CL}
\tag{(CL)}
\rho(y)-\rho(y+s)=\phi_s(y+1)-\phi_s(y)+c_s\ .
\end{equation}

The group $\CG(X)$ associated to the system is 
the group of transformations of $X=Y\times \T^d$ of the form
$$
 (y,h)\mapsto (y+s,h+\phi_s(y))
$$
where $s$ and $\phi_s$ satisfy~\ref{eq:CL}.

\subsection{Structure theorem in general}

We give a short outline of the steps needed to 
complete the proof of the Structure Theorem for $k\geq 3$.
We have that $Y=Z_{k-1}(X)$ is a system of order $k-1$, 
$X=Y\times \T^d$, $T(y,h)=(Ty,h+\rho(y))$,  and 
$\rho\colon Y\to\T^d$ satisfies the functional equation~\eqref{eq:funceq}.
By the induction hypothesis $Y=\limproj Y_i$ where 
each $Y_i=G_i/\Gamma_i$ is a $(k-1)$-step nilsystem.

We first show that the cocycle $\rho$ is {\em cohomologous} 
to a cocycle measurable with respect to $\CY_i$ for some $i$, 
meaning that the difference between the two cocycles is a 
coboundary. 
This reduces us to the case that $\rho$ is measurable with respect to 
some $\CY_i$, and so we
can assume that $Y=\CY_i$ for some $i$.  Thus 
$Y$ is a $(k-1)$-step nilsystem and 
we can assume that $Y=G/\Gamma$ with $G=\CG(Y)$.

We then use the functional equation to lift every 
transformation 
$S\in G$ to a transformation of $X$ belonging to $\CG(X)$. 
Starting with the case $S\in 
G_{k-1}$, we move up the lower central series of $G$.
Lastly we show that 
we obtain sufficiently many elements of the group $\CG(X)$ 
in this way.

\subsection{Relations to the finite case}

The seminorms $\nnorm\cdot_k$ play the same role that the Gowers
norms play in Gowers's proof~\cite{gowers} of Szemer\'edi's Theorem and 
in Green and Tao's proof~\cite{GT} that the primes 
contain arbitrarily long arithmetic progressions.  
We let  $U_k$ denote the $k$-th Gowers norm.  
For the finite system $\Z/N\Z$, $\nnorm f_k=\norm f_{U_k}$.
Furthermore, $\norm\cdot_{U_k}$ is a norm, not only a seminorm.
The analog of Lemma~\ref{lemma:charseminorm} 
is that if $\|f_0\|_\infty, 
\|f_1\|_\infty, \dots, \|f_k\|_\infty\leq 1$, then there exists 
some constant $C_k > 0$ such that 
$$
 \bigl|\E\bigl(f_0(x)f_1(x+y)\dots f_k(x+ky)\mid x,y\in\Z/p\Z\bigr)
\bigr|
 \leq
C_k\min_{0\leq j\leq k}\norm {f_j}_{U_k} \ .
$$

Other parts of the program are not as easy to translate 
to the finite setting.  Consider defining a factor of 
the system using the seminorms.    
If $p$ is prime, then $\Z/p\Z$ has no nontrivial factor 
and so there is no factor of $\Z/p\Z$ playing the 
role of the factor $Z_k$, meaning there is no factor with
$$
\E(f\mid \CZ_k)=0 \text{ if and only if } \norm f_{U_k}=0\ .
$$
Instead, the corresponding results have a different flavor: if $\norm
f_{U_k}$ is large in some sense, then $f$ has large conditional expectation
on some (noninvariant) $\sigma$-algebra or it has large correlation
with a function of some particular class.
Although we have a complete characterization of the seminorms 
$\nnorm \cdot_k$ (and so the factors 
$Z_k$)  in terms of nilmanifolds, there are only partial combinatorial 
characterizations in this direction 
(see~\cite{GT2},~\cite{GT3} and~\cite{GT4}).

\section{Other patterns}
\label{sec:other}

\subsection{Commuting transformations}
\label{sec:commuting}

Ergodic theory has been used to detect other patterns that 
occur in sets of positive upper density, using Furstenberg's 
Correspondence Principle and an appropriately 
chosen strengthening of Furstenberg 
multiple recurrence.  A first example is for 
commuting transformations:

\begin{theorem}[Furstenberg and Katznelson~\cite{FuK1}]
Let $(X,\CX,\mu)$ be a probability measure space, let 
$k\geq 1$ be an integer,  and assume that $T_j\colon X\to X$ 
are commuting measure preserving transformations for $j=1, 2, \ldots, k$, 
then for all $A\in\CX$ with $\mu(A)>0$, there exist infinitely many 
$n\in\N$ such that 
\begin{equation}
\label{eq:commuting}
\mu(A\cap T_1^{-n}A\cap T_2^{-n}A\cap\ldots\cap T_k^{-n}A)> 0 \ .
\end{equation}
\end{theorem}

(In~\cite{FuK2}, Furstenberg and Katznelson proved a 
strengthening of this result, showing that one can place 
some restrictions on the choice of $n$; we do not discuss 
these ``IP'' versions of the theorems given in the sequel.)
Via correspondence, a multidimensional version of 
Szemer\'edi's Theorem follows: if $E\subset\Z^r$ has positive 
upper density and $F\subset\Z^r$ is a finite subset, then 
there exist $z\in \Z^r$ and $n\in\N$ such that $z+nF\subset E$.

Again, this theorem is proven by showing that the 
associated $\liminf$ of the average of 
the quantity in Equation~\eqref{eq:commuting} 
is positive.  Again, it is natural to ask 
if the limit
$$
\lim_{N\to\infty}\frac{1}{N}\sum_{n=0}^{N-1}
\mu(A\cap T_1^{-n}A\cap \ldots\cap T_k^{-n}A) 
$$ 
exists in $L^2(\mu)$ for commuting maps $T_1, \ldots, T_k$.
Only partial results are known.  For $k=2$, 
Conze and Lesigne (\cite{CL3},~\cite{CL4}) proved convergence.  For 
$k\geq 3$, the only known results rely on strong 
hypotheses of ergodicity:  

\begin{theorem}[Frantzikinakis and Kra~\cite{FK2}]
\label{th:FKcommuting}
Let $k\in\N$ and assume that $T_1, T_2, \ldots, T_k$
are commuting invertible ergodic measure preserving
transformations of a measure space $(X, \CX, \mu)$ such that
$T_iT_j^{-1}$ is ergodic for all $i,j\in \{1, 2, \ldots, k\}$ with
$i\neq j$.  If  $f_1, f_2, \ldots, f_k\in L^\infty(\mu)$ the
averages,
$$
\frac{1}{N}\sum_{n=0}^{N-1}T_1^nf_1\cdot
T_2^nf_2\cdot\ldots \cdot T_k^nf_k
$$
converge in $L^2(\mu)$ as $N\to\infty$.
\end{theorem}

The idea is to prove an analog of Lemma~\ref{lemma:charseminorm} 
for commuting transformations, thus reducing the problem 
to working in a nilsystem.  The factors $Z_k$ that are characteristic 
for averages along arithmetic progressions are also characteristic 
for these particular averages of commuting transformations.  
Without the strong hypotheses of ergodicity, this no longer 
holds and the general case remains open.  

\subsection{Averages along cubes}
Another type of average is along $k$-dimensional cubes, 
the natural objects that arise in the definition of the 
seminorms.  For example, a $2$-dimensional cube is an 
expression of the form:
$$f(x)f(T^mx)f(T^nx)f(T^{m+n}x)\ .
$$

In~\cite{Berg5}, Bergelson 
showed  the existence in $L^2(\mu)$ of 
$$
\lim_{N\to\infty}\frac{1}{N^2}\sum_{n,m=0}^{N-1}
T^nf_1\cdot T^mf_2\cdot T^{n+m}f_3 \ ,
$$
where $f_1, f_2, f_3\in L^\infty(\mu)$.  
Similarly, one can define a $3$-dimensional cube:
$$
f_1(T^mx)f_2(T^nx)f_3(T^{m+n}x)f_4(T^px)f_5(T^{m+p}x)f_6(T^{n+p}x)f_7(T^{m+n+p}x)
$$
and existence of the limit of the average of this 
expression $L^2(\mu)$ for bounded functions 
$f_1, f_2, \ldots, f_7$ was shown in~\cite{HK3}.  

More generally, this theorem holds for cubes of $2^{k}-1$ functions.
Recalling the notation of Section~\ref{sec:charfac}, we have 
for $\epsilon =\epsilon _1\dots \epsilon _k\in \{0,1\}^k$ and
$\bn=(n_1,\dots,n_k)\in\Z^k$,
$$
\epsilon \cdot \bn=\epsilon _1n_1+\epsilon _2n_2+\dots+\epsilon _kn_k\ ,
$$
and $\bzero$ denotes the element $00\dots 0$ of $\{0,1\}^k$.  
We have:
\begin{theorem}[Host and Kra~\cite{HK4}]
\label{th:cubes}
Let $(X, \CX, \mu, T)$ be a system, let $k\geq 1$ be an integer, 
and let $f_\epsilon $, $\epsilon \in\{0,1\}^k\setminus\{\bzero\}$,
be $2^k-1$ bounded functions on $X$.
Then the averages
$$
\frac {1}{N^k}\ \cdot \
\sum_{\bn\in[0,N-1]^k}\
\prod_{\substack{\epsilon \in\{0,1\}^k\\ \epsilon\neq\bzero}}
T^{\epsilon \cdot \bn}f_\epsilon 
$$
converge in $L^2(\mu)$ as $N\to\infty$. 
\end{theorem}

The same result holds for 
translated averages, meaning the average for 
$\bn\in[M_1,N_1]\times\dots\times [M_k,N_k]$, as $N_1-M_1$, \dots,
$N_k-M_k\to\infty$.

By Furstenberg's Correspondence Principle, this translates 
to a combinatorial statement.  
A subset $E\subset\Z$ is {\em syndetic} if $\Z$ can be covered 
by finitely many translates of $E$.  In other words, there exists 
$N> 0$ such that every interval of size $N$ contains at least one 
element of $E$.  (Thus it is natural to refer to a 
syndetic set in the integers as a set with {\em bounded gaps}.)
More generally, 
$E\subset \Z^k$ is {\em syndetic} if there exists an integer 
$N>0$ such that 
$$E\cap \bigl([M_1, M_1+N]\times \ldots\times 
[M_k, M_k+N]\bigr)\neq\emptyset
$$ 
for all $M_1, \ldots, M_k\in \Z$.

Restricting Theorem~\ref{th:cubes} 
to indicator functions, the limit of the averages
$$
\prod_{i=1}^k\frac{1}{\scriptstyle N_i-M_i}\ \cdot \
\sum_{\substack{n_1\in[M_1,N_1], \ldots, n_k\in[M_k,N_k]}}\
\mu\bigl(\bigcap_{\epsilon \in \{0,1\}^k}T^{\epsilon \cdot \bn}A\bigr)
$$
exists and is greater than or equal to $\mu(A)^{2^k}$ when
$N_1- M_1, \ldots, N_k- M_k\to\infty$.
Thus for every $\varepsilon>0$,
$$
  \Bigl\{
  \bn\in\Z^k:
  \mu\bigl(\bigcap_{\epsilon \in \{0,1\}^k}T^{\epsilon \cdot \bn}
  A\bigr)> \mu(A)^{2^k}-\varepsilon
  \Bigr\}
$$
of $\Z^k$ is syndetic.

By the Correspondence Principle, we have that if 
$E\subset\Z$ has upper density 
$d^*(E) > \delta > 0$ and $k\in\N$, then 
$$
\Bigl\{\bn\in \Z^k\colon
d^*\bigl(\bigcap_{\epsilon\in \{0,1\}^k}(E+\epsilon\cdot\bn)\bigr)
\geq \delta^{2^k}\Big\}
$$ 
is syndetic.  

\subsection{Polynomial patterns}

In a different direction, one can restrict the iterates 
arising in Furstenberg's multiple recurrence.  
A natural choice is polynomial iterates, and the 
corresponding combinatorial statement is that 
a set of integers with positive upper density contains 
elements who differ by a polynomial:

\begin{theorem}[S\'ark\"ozy~\cite{sarkozy}, Furstenberg~\cite{Fbook}]
If $E\subset\N$ has positive upper density and $p\colon\Z\to\Z$ is a
polynomial with $p(0)=0$, then there
exist $x,y\in E$  and $n\in \N$ such that
$x-y=p(n)$.
\end{theorem}

As for arithmetic progressions, 
Furstenberg's proof relies on an averaging theorem:

\begin{theorem}[Furstenberg~\cite{Fbook}]
Let $(X,\CX,\mu,T)$ be a system, let $A\in\CX$ with 
$\mu(A)>0$ and let $p\colon\Z\to\Z$ be a polynomial with $p(0)=0$. Then
$$ 
\liminf_{N\to\infty}\frac{1}{N}\sum_{n=0}^{N-1}
    \mu(A\cap T^{-p(n)}A)>0 \ .
$$
\end{theorem}

The multiple polynomial recurrence theorem, simultaneously generalizing 
this single polynomial result and Furstenberg's multiple recurrence, 
was proven by Bergelson and Leibman:
\begin{theorem}[Bergelson and Leibman~\cite{BL}]
\label{th:BLpoly}
Let $(X,\CX,\mu,T)$ be a system,  
let $A\in\CX$ with $\mu(A)>0$, and let $k\in\N$.  
If $p_1, p_2, \ldots, p_k\colon\Z\to\Z$ are polynomials 
with $p_j(0)=0$ for 
$j=1, \ldots, k$, then 
\begin{equation}
\label{eq:polyav}
\liminf_{N\to\infty}\frac 1N\sum_{n=0}^{N-1}\mu\bigl(A\cap T^{-p_1(n)}A 
\cap\dots\cap 
 T^{-p_k(n)}A\bigr)  > 0 \ .
\end{equation}
\end{theorem}
By the Correspondence Principle, 
one immediately deduces a polynomial Szemer\'edi Theorem: 
if $E\subset\Z$ has positive upper density, then it contains 
arbitrary polynomial patterns, meaning there exists $n\in\N$ such that 
$$
x, x + p_1(n), x+p_2(n), \ldots, x+p_k(n)\in E \ .
$$
(More generally, Bergelson and Leibman proved a version of 
Theorem~\ref{th:BLpoly} 
for commuting transformations, with a multidimensional 
polynomial Szemer\'edi Theorem as a corollary.)

Again, it is natural to ask if the $\liminf$ 
in~\eqref{eq:polyav} is actually a limit.  
A first result in this direction 
was given by Furstenberg and Weiss~\cite{FW}, who proved convergence  
in $L^2(\mu)$ of
$$
 \frac 1N\sum_{n=0}^{N-1} T^{n^2}f_1\cdot T^nf_2
$$
and
$$
 \frac 1N\sum_{n=0}^{N-1} T^{n^2}f_1\cdot T^{n^2+n}f_2
$$
for bounded functions $f_1, f_2$.  

The proof of convergence for general polynomial 
averages uses the technology 
of the seminorms, reducing to the same 
characteristic factors $Z_k$ that can be described using 
nilsystems, as for averages along arithmetic progressions:
\begin{theorem}[Host and Kra~\cite{HK5}, Leibman~\cite{Leibman1}]
Let $(X,\CX,\mu,T)$ be a system, $k\in\N$, and 
$f_1,f_2,\dots,f_k\in L^\infty(\mu)$. 
Then for any polynomials 
$p_1, p_2, \dots, p_k\colon\Z\to\Z$, the averages
\begin{equation*}
\frac 1N\sum_{n=0}^{N-1}
T^{p_1(n)}f_1\cdot T^{p_2(n)}f_2\cdot\ldots\cdot T^{p_k(n)}f_k
\end{equation*}
converge in $L^2(\mu)$.
\end{theorem}

Recently, Johnson~\cite{johnson} has shown that 
under similar strong ergodicity conditions to those 
in Theorem~\ref{th:FKcommuting}, one can generalize 
this and prove $L^2(\mu)$-convergence of the polynomial averages 
for commuting transformations:
\begin{equation*}
\frac 1N\sum_{n=0}^{N-1}
T_1^{p_1(n)}f_1\cdot T_2^{p_2(n)}f_2\cdot\ldots\cdot T_k^{p_k(n)}f_k
\end{equation*}
for $f_1, f_2, \ldots, f_k\in L^\infty(\mu)$.

For a {\em totally ergodic system} (meaning that 
$T^n$ is ergodic for all $n\in\N$), 
Furstenberg and Weiss showed a 
stronger result, giving an explicit and simple formula for 
the limit:
$$
\frac{1}{N}\sum_{n=0}^{N-1} T^nf_1\cdot T^{n^2}f_2
 \to \int f_1 \, d\mu \cdot \int f_2 \, d\mu
$$ 
in $L^2(\mu)$.  

Bergelson~\cite{Berg4} asked whether the same result holds for $k$
polynomials of different degrees, meaning that the limit of 
the polynomial average for a totally ergodic system is 
the product integrals.  We show that the answer is yes under 
a more general condition.
A family of polynomials 
$p_1, p_2, \ldots, p_k\colon\Z\to\Z$
is {\em rationally independent} if for all integers $m_1,
\ldots, m_{k}$ with at least some $m_{j}\neq 0 $,
the polynomial $\sum_{j=1}^k m_jp_j(n)$ is not constant.
We show:

\begin{theorem}[Frantzikinakis and Kra~\cite{FK1}]
Let $(X, \CX, \mu, T)$ be a totally ergodic 
system, let $k\geq 1$ be an integer, and assume that $p_1,p_2
\ldots, p_k\colon\Z\to\Z$ are rationally independent 
polynomials. If $f_1, f_2, \ldots, f_k \in L^\infty(\mu)$,
$$
\lim_{N\to\infty}\Big\Vert
\frac{1}{N}\sum_{n=0}^{N-1} T^{p_1(n)}f_1\cdot T^{p_2(n)}\cdot\ldots\cdot
T^{p_k(n)}f_k 
- \prod_{i=1}^k\int f_i\,d\mu
\Big\Vert_{L^2(\mu)} = 0 \ .
$$
\end{theorem}

As a corollary, if $(X,\CX,\mu,T)$ is  totally ergodic, 
$\{1,p_1,\ldots, p_k \}$ are rationally independent 
polynomials taking on integer values on the integers, and 
$A_0, A_1, \ldots, A_k\in\CX$ with $\mu(A_i)>0$, $i=0,\ldots,k$,
then
$$
\mu(A_0\cap T^{-p_1(n)}A_1\cap \ldots \cap
 T^{-p_k(n)}A_k)>0
$$
for some $n\in\N$.
Thus in a totally ergodic system, one can strengthen 
Bergelson and Leibman's multiple polynomial recurrence 
theorem, allowing the sets $A_i$ to be distinct, and 
allowing the polynomials $p_i$ to have nonzero constant term.  
It is not clear if this has a combinatorial interpretation.  

\section{Strengthening Poincar\'e recurrence}
\subsection{Khintchine recurrence}

Poincar\'e recurrence states that a set of positive measure 
returns to intersect itself infinitely often.  One way 
to strengthen this is to ask that the set return to itself 
often with `large' intersection.  
Khintchine made this notion precise, showing that 
large self intersection occurs on a syndetic set:

\begin{theorem}[Khintchine~\cite{Khintchine}]
Let $(X, \CX, \mu, T)$ be a 
system, let $A\in\CX$ have $\mu(A) > 0$, and let $\varepsilon 
> 0$.  Then
$$\{n\in\Z: \mu(A\cap T^nA)> \mu(A)^2-\varepsilon\}$$
is syndetic.   
\end{theorem}

It is natural to ask for a 
simultaneous generalization of Furstenberg 
Multiple Recurrence and Khintchine Recurrence.  More 
precisely, if $(X, \CX,\mu, T)$ is a  system, $A\in\CX$ has positive 
measure, $k\in\N$, and $\varepsilon>0$, is the set
$$
  \bigl\{n\in\Z\colon\mu\bigl(A\cap T^nA\cap\dots\cap
  T^{kn}A\bigr)>\mu(A)^{k+1}-\varepsilon\bigr\}
$$
syndetic?

Furstenberg Multiple Recurrence implies that there 
exists some constant $c = c(\mu(A))>0$ such that 
$$
\{n\in\Z\colon \mu(A\cap T^nA\cap \ldots\cap T^{kn}A)> 
c \}
$$ 
is syndetic.  But to generalize Khintchine 
Recurrence, one needs $c=\mu(A)^{k+1}$.  
It turns out that the answer depends on the length $k$ of the arithmetic progression.

\begin{theorem}[Bergelson, Host and Kra~\cite{BHK}]
\label{th:BHK}
Let $(X,\CX,\mu,T)$ be an ergodic system and let $A\in\CX$.
Then for every $\varepsilon>0$, the sets
$$
    \bigl\{n\in\Z : \mu(A\cap T^nA\cap
      T^{2n}A)>\mu(A)^3-\varepsilon\bigr\}
$$
and
$$
          \bigl\{n\in\Z : \mu(A\cap T^nA\cap T^{2n}A\cap T^{3n}A)
>\mu(A)^4-\varepsilon\bigr\}
$$
are syndetic.
\end{theorem}
Furthermore, this result fails on average, meaning that the
average of the left hand side expressions is not greater than 
$\mu(A)^3-\varepsilon$ or $\mu(A)^4-\varepsilon$, respectively.  

On the other hand, based on an example of Ruzsa 
contained in the appendix of~\cite{BHK}, we have:
\begin{theorem}[Bergelson, Host and Kra~\cite{BHK}]
\label{th:BHK2}
There exists an ergodic system $(X,\CX,\mu,T)$ and for all 
$\ell\in\N$ there exists a  set $A=A(\ell)\in\CX$
with $\mu(A)>0$ such that
$$
\mu(A\cap T^nA\cap T^{2n}A\cap T^{{3n}}A
\cap T^{4n}A)\leq \mu(A)^{\ell}/2
$$
for every integer $n\neq 0$.
\end{theorem}

We now briefly outline the major ingredients in 
the proofs of these theorems.  

\subsection{Positive ergodic results}

We start with the ergodic results needed to prove Theorem~\ref{th:BHK}.  
Fix an integer $k\geq 1$, an ergodic  system $(X,\CX, \mu,
T)$, and $A\in\CX$ with $\mu(A)>0$. 
The key ingredient 
is the study of the {\em multicorrelation sequence}
$$
\mu\bigl(A\cap T^n A\cap T^{2n} A \cap\ldots\cap T^{kn}A\bigr) \ .
$$
More generally, for a real valued 
function $f\in L^\infty(\mu)$, we consider the {\em 
multicorrelation sequence}
$$
I_{f}(k, n) := \int f\cdot T^nf\cdot T^{2n}f\cdot\ldots\cdot
T^{kn}f\,d\mu(x) \ .
$$

When $k =1$, Herglotz's Theorem implies that the correlation
sequence $I_{f}(1, n)$ is the
Fourier transform of some positive measure $\sigma=\sigma_f$ on the torus
$\T$:
$$
I_{f}(1, n) = \wh{\sigma}(n):= \int_{\T} e^{2\pi i nt}\,
d\sigma(t) \ .
$$
Decomposing the measure $\sigma$ into its
continuous part $\sigma^c$ and its discrete part $\sigma^d$, 
can write the multicorrelation sequence $I_{f}(1, n)$ as the sum
of two sequences
$$
I_{f}(1, n) = \wh{\sigma^c}(n) + \wh{\sigma^d}(n) \ .
$$
The sequence $\{\wh{\sigma^c}(n)\}$ {\em tends to $0$ in density}, 
meaning that 
$$
\lim_{N\to\infty} \sup_{M\in\Z}
     \frac 1M\sum_{n=M}^{M+N-1}\lvert \wh{\sigma^c(n)}\rvert =0\ .
$$
Equivalently, for any $\varepsilon>0$, the 
upper Banach density\footnote{The {\em upper Banach 
density} $\overline{d}(E)$ of a set $E\subset\Z$ is defined by 
$\overline{d}(e) = \lim_{N\to\infty}\sum_{M\in\Z}\frac{1}{N}
\vert E\cap[M, M+N-1]\vert$.} of the set
$\{n\in\Z: \lvert \wh{\sigma^c(n)}\rvert>\varepsilon\}$
is zero.  
The sequence $\{\wh{\sigma^d}(n)\}$ is {\em almost periodic},
meaning that there exists a compact abelian group $G$, a continuous
real valued function $\phi$ on $G$, and $a\in G$ such that $\wh{\sigma^d}(n) =
\phi(a^n)$ for all $n$.

A compact abelian group can be approximated by  a compact abelian Lie
group.  Thus any almost periodic sequence can be uniformly
approximated by an almost periodic sequence arising from a compact
abelian Lie group.

In general, however, for higher $k$ the answer is more complicated.
We find a similar decomposition for the multicorrelation sequences
$I_f(k,n)$ for $k\geq 2$.
The notion of an almost periodic sequence is replaced by that of a
{\em nilsequence}:  
for an integer $k\geq 2$,  
a $k$-step nilmanifold $X=G/\Gamma$, 
a continuous real (or complex) valued function $\phi$ on 
$G$, $a\in G$,  and $e\in X$, the 
sequence $\{\phi(a^n\cdot e)\}$ is called a {\em basic $k$-step
nilsequence}.  
A {\em $k$-step nilsequence} is a uniform limit of basic
$k$-step nilsequences.

It follows that a $1$-step nilsequence is the same as an almost periodic
sequence.  An inverse limit of compact abelian Lie groups is a compact group. 
However an inverse limit of $k$-step nilmanifolds is not, in general, the  
homogeneous space of some locally compact group, and so for higher $k$, 
the decomposition result must take into account the 
uniform limits of basic nilsequences.  We have: 
\begin{theorem}[Bergelson, Host and Kra~\cite{BHK}]
Let $(X,\CX,\mu,T)$ be an ergodic system, $f\in L^\infty(\mu)$ and
$k\geq 1$ an integer. The sequence $\{I_f(k,n)\}$ is the sum of
a sequence tending to zero in density and a $k$-step
nilsequence.
\end{theorem}

Finally, we explain how this result can be used to prove 
Theorem~\ref{th:BHK}.  
Let $\{a_n\}_{n\in\Z}$ be a bounded sequence of real numbers.
The {\em syndetic supremum} of this sequence is defined to be
$$
\sup
\Bigl\{ c\in \R\colon \{n\in\Z\colon a_n > c \} \text{ is syndetic }\Bigr\}\ .
$$
Every nilsequence $\{a_n\}$ is uniformly recurrent. In particular, if
$S=\sup(a_n)$ and $\varepsilon>0$, then
$\{n\in\Z\colon a_n\geq S-\varepsilon\}$ is syndetic.

If $\{a_n\}$ and $\{b_n\}$ are two sequences of real numbers such that 
$a_n-b_n$ tends to $0$ in density, then the two sequences have the 
same syndetic supremum.
Therefore the syndetic supremums of the sequences
$$\{\mu(A\cap T^nA\cap T^{2n}A)\}$$ and
$$\{\mu(A\cap T^nA\cap T^{2n}A\cap T^{3n}A)\}$$
are equal to the supremum of the associated nilsequences, and we are
reduced  to showing that they are
greater than or equal to $\mu(A)^3$ and $\mu(A)^4$, respectively.

\subsection{Nonergodic counterexample}

Ergodicity is not needed for Khintchine's Theorem, but 
is essential for Theorem~\ref{th:BHK}:
\begin{theorem}[Bergelson, Host, and Kra~\cite{BHK}]
\label{th:counter}
There exists a (nonergodic) system $(X,\CX,\mu,T)$, and for every
$\ell\in\N$ there exists $A\in\CX$  with $\mu(A) > 0$ such that
$$
    \mu\bigl(A\cap T^nA\cap T^{2n}A\bigr)\leq \frac 12\mu(A)^\ell\ .
$$
for integer $n\neq 0$.
\end{theorem}

Actually there exists 
a set $A$ of arbitrarily small positive measure with
$$\mu\bigl(A\cap T^nA\cap T^{2n}A\bigr)\leq \mu(A)^{-c\log(\mu(A))}$$
for every integer $n\neq 0$ and for 
some positive universal constant $c$.

The proof is based on Behrend's construction of a set containing 
no arithmetic progression of length $3$:
\begin{theorem}[Behrend~\cite{behrend}]
For all $L\in\N$, there
exists a subset $E\subset\{0,1,\dots,L-1\}$ 
having more than $L\exp(-c\sqrt{\log L})$ elements
that does not contain any nontrivial arithmetic progression of length
$3$.
\end{theorem}

\begin{proof}(of Theorem~\ref{th:counter})
Let $X=\T\times \T$, with Haar measure
$\mu=m\times m$ and transformation $T\colon X\to X$ given by
$T(x,y)=(x,y+x)$.

Let $E\subset\{0,1,\dots,$ $L-1\}$, not
containing any nontrivial arithmetic progression of length
$3$. Define
$$
    B=\bigcup_{j\in E}\bigl[\frac j{2L},\frac j{2L}+\frac 1{4L}\bigr)\ ,
$$
which we consider as a subset of the torus and $A=\T\times B$.

For every integer $n\neq 0$, we have $T^n(x,y)=(x,y+nx)$ and
\begin{align*}
\mu\bigl(A\cap T^nA\cap T^{2n}A\bigr) &
= \iint_{\T\times\T}\one_B(y)\one_B(y+nx)\one_B(y+2nx)\,dm(y)\,dm(x) \\
& = \iint_{\T\times\T}\one_B(y)\one_B(y+x)\one_B(y+2x)\,dm(y)\,dm(x)\ .
\end{align*}
Bounding this integral, we have that:
\begin{align*}
\mu\bigl(A\cap T^nA\cap T^{2n}A\bigr) &  = 
\iint_{\T\times\T}\one_B(y)\one_B(y+x)\one_B(y+2x)\,dm(x)\,dm(y)\\
& \leq \frac{m(B)}{4L}\ .
\end{align*}

By Behrend's Theorem, we can choose the set $E$ with cardinality on the order
of $L\exp(-c\sqrt{\log L})$. Choosing $L$ sufficiently large,
a simple computation gives the statement.  
\end{proof}

For longer arithmetic progressions, the counterexample 
of Theorem~\ref{th:BHK2} is 
based on a construction of Ruzsa.  
When $P$ is a nonconstant integer polynomial of degree $\leq 2$,
the subset
$$
    \bigl\{P(0),P(1),P(2),P(3),P(4)\bigr\}
$$
of $\Z$ is called a {\em quadratic configuration of $5$ terms},
written QC5 for short.

Any QC5 contains at least $3$ distinct elements.
An arithmetic progression of length $5$ is a QC5,
corresponding to a polynomial of degree $1$.

\begin{theorem}[Ruzsa~\cite{BHK}]
For all $L\in\N$, there exists  a subset $E\subset\{0,1,\dots,$ $L-1\}$
having more than $L\exp(-c\sqrt{\log L})$ elements that does not
contain  any QC5.
\end{theorem}

Based on this, we show:
\begin{theorem}[Bergelson, Host and Kra~\cite{BHK}]
There exists an ergodic system $(X,\CX, \mu,T)$ and, for every 
$\ell\in\N$, there exists $A\in\CX$ with $\mu(A) > 0$ such that
$$
\mu(A\cap T^nA\cap T^{2n}A\cap T^{3n}A\cap T^{4n}A)
\leq \frac 12\mu(A)^\ell
$$
for every integer $n\neq 0$.
\end{theorem}

Once again, proof gives the estimate 
$\mu(A)^{-c\log(\mu(A))}$, for some constant $c> 0$.  

The construction again involves a simple example: 
$\T$ is the torus with Haar measure $m$,   
$X=\T\times \T$, and  $\mu=m\times m$. Let $\alpha\in\T$ be 
irrational and let $T\colon X\to X$ be
$$
 T(x,y)=(x+\alpha,y+2x+\alpha)\ .
$$

Combinatorially this example becomes:
for all $k\in\N$, there exists $\delta>0$ such that for
infinitely many integers $N$, there is a subset $A\subset
\{1,\dots, N\}$ with $\lvert A\rvert\geq\delta N$ that contains no
more than
$\frac 12 \delta^kN$ arithmetic progressions of length $\geq 5$ 
with the same difference.

\subsection{Combinatorial consequences}
Via a slight modification of the Correspondence Principle, 
each of these results translates to a combinatorial statement.  
For $\varepsilon > 0$ and $E\subset\Z$ with positive upper Banach 
density, consider the set 
\begin{equation}
\label{eq:intersection}
\{n\in\Z\colon \overline{d}(E\cap (E+n)\cap (E+2n)\cap \ldots\cap (E+kn)) 
\geq \overline{d}(E^{k+1})-\varepsilon\} \ .
\end{equation}
From Theorems~\ref{th:BHK} and~\ref{th:BHK2}, for $k=2$ and for $k=3$, 
this set is syndetic, while for $k\geq 4$ there exists a set of integers 
$E$ with positive upper Banach density such that the set 
in~\eqref{eq:intersection} is empty.  

We can refine this a bit further.  
Recall the notation from Szemer\'edi's Theorem: 
for every $\delta > 0$ and $k\in\N$, there
exists $N(\delta, k)$ such that for all $N> N(\delta, k)$,
every subset of $\{1,\dots, N\}$ with at least $\delta N$
elements contains an arithmetic progression of length $k$.

For an arithmetic progression $\{a, a+s, \dots,
a+(k -1)s\}$, $s$ is the {\em difference} of the progression.
Write $\lfloor x\rfloor$ for integer part of $x$.
From Szemer\'edi's Theorem, we can 
deduce that every subset $E$ of $\{1,\ldots, N\}$
with at least $\delta N$ elements contains at least
$\lfloor cN^2\rfloor$
arithmetic progressions of length $k$, where
$c = c(k,\delta)> 0$ is a constant. 
Therefore the set $E$
contains at least $\lfloor c(k,\delta)N\rfloor$ progressions of 
length $k$ with
the same difference.

The ergodic results of Theorem~\ref{th:BHK} give 
some improvement for $k=3$ and $k=4$ (see~\cite{BHK} for the precise 
statement).  
For $k=3$, this was strengthened by Green:

\begin{theorem}[Green~\cite{green}]
For all $\delta, \varepsilon >0$, 
there exists $N_0(\delta, \varepsilon)$ such that for all $N>N_0(\delta, 
\varepsilon)$ and any $E\subset\{1,\ldots, N\}$ with $|E|\geq \delta N$, 
$E$ contains at least $(1-\varepsilon)\delta^3N$ arithmetic progressions 
of length $3$ with the same difference.  
\end{theorem}

On the other hand, the similar bound for longer progressions 
with length $k\geq 5$ does not hold. 
The proof in~\cite{BHK}, based on an example of Rusza,  
does not use ergodic theory.  We show that for 
all $k\in\N$, there exists $\delta>0$ such that for
infinitely many $N$, there exists a subset $E$ of
$\{1,\dots, N\}$ with $\lvert E\rvert\geq\delta N$ that contains no
more than $\frac 12 \delta^kN$ arithmetic progressions of length $\geq 5$ 
with the same step.

\subsection{Polynomial averages}

One can ask if similar lower bounds hold for the polynomial 
averages.  For independent polynomials, using the fact that 
the characteristic factor is the Kronecker factor, we 
can show: 
\begin{theorem}[Frantzikinakis and Kra~\cite{FK3}]
Let $k\in\N$, $(X,\CX,\mu,T)$ be a system, $A\in\CX$, 
and let $p_1, p_2, \ldots, p_k\colon\Z\to\Z$ be rationally independent 
polynomials with $p_i(0)=0$ for $i=1, 2m \ldots, k$.
Then for every $\varepsilon>0$, the set
$$
      \bigl\{n\in\Z : \mu(A\cap T^{p_1(n)}A\cap T^{p_2(n)}\cap\ldots\cap
      T^{p_k(n)}A) 
>\mu(A)^{k+1}-\varepsilon\bigr\}
$$
is syndetic.
\end{theorem}
Once again, this result fails on average.

Via Correspondence, analogous to the results of~\eqref{eq:intersection}, 
we have that 
for $E\subset\Z$ and 
rationally independent polynomials $p_1, p_2, \ldots, p_k\colon\Z\to\Z$ 
with $p_i(0)=0$ 
for $i=1, 2, \ldots, k$, then for all $\varepsilon > 0$, the set 
$$
\{n\in\Z\colon \overline{d}
\bigl(E\cap (E+p_1(n))\cap\ldots \cap(E+p_k(n))\bigr) 
\geq \overline{d}(E)^{k+1}-\varepsilon\}
$$
is syndetic.

Moreover, in~\cite{FK3} we strengthen this and show that 
there are many configurations with the same $n$ giving 
the differences:
if $p_1, p_2, \ldots, p_k\colon\Z\to\Z$ 
are rationally independent polynomials with $p_i(0)=0$ 
for $i=1, 2, \ldots, k$, then for all $\delta, \varepsilon > 0$, 
there exists $N(\delta, \varepsilon)$ such that for all $N> N(\delta, 
\varepsilon)$ and any subset $E\subset\{1, \ldots, N\}$ 
with $|E|\geq \delta N$ contains at least 
$(1-\varepsilon)\delta^{k+1}N$ configurations of the form 
$$\{x, x+p_1(n), x+p_2(n), \ldots, x+p_k(n)\}$$ for a fixed $n\in\N$.

\bibliographystyle{amsalpha}

\begin{thebibliography}{9999}

\bibitem{behrend}
F.~A.~Behrend.  On sets of integers which contain no three 
in arithmetic progression.  {\em Proc. Nat. Acad. Sci.} {\bf 23}
(1946), 331--332.  

\bibitem{Berg3}
V.~Bergelson.
Weakly mixing PET.
{\em Erg. Th. \& Dyn. Sys.} {\bf 7}  (1987), 337--349.

\bibitem{Berg4}
V.~Bergelson.  Ergodic Ramsey theory an update. {\em Ergodic
Theory of ${\Z}^d$-actions} (Eds.: M. Pollicott, K. Schmidt).
Cambridge University Press, Cambridge (1996), 1--61.

\bibitem{Berg5}
V.~Bergelson.  
The multifarious Poincaré recurrence theorem.
{\em Descriptive set theory and dynamical systems} 
(Marseille-Luminy, 1996), 
Cambridge University Press, Cambridge (2000), 31--57.

\bibitem{BHK}
V.~Bergelson, B.~Host and B.~Kra, with an Appendix
by I.~Ruzsa.  Multiple recurrence
and nilsequences.  {\it Inventiones Math.} 160 (2005), 261-303.

\bibitem{BL}
V.~Bergelson and A.~Leibman.
Polynomial extensions of van der Waerden's and Szemer\'edi's
theorems.  {\em J. Amer. Math. Soc.} {\bf 9} (1996), 725--753.

\bibitem{Bourgain2}
J.~Bourgain.
On the maximal ergodic theorem for certain subsets of the
positive integers.  {\em Isr. J. Math.} {\bf 61} (1988), 39--72.


\bibitem{CL3}
J.-P. Conze and E.~Lesigne.
Sur un th\'{e}or\`{e}me ergodique pour des mesures diagonales.
{\em Publications de l'Institut de Recherche de Math\'ematiques de
  Rennes, Probabilit\'es} 1987.

\bibitem{CL4}
J.-P. Conze and E.~Lesigne.
Sur un th\'{e}or\`{e}me ergodique pour des mesures diagonales.
{\em C. R. Acad. Sci. Paris S\'{e}rie {I}} {\bf 306}  (1988), 491--493.

\bibitem{CFS}
I.~Cornfeld, S.~Fomin and Ya.~Sinai. {\em Ergodic Theory}.  
Springer-Verlag, Berlin, Heidelberg, New York 1982.


\bibitem{ET}
P.~Erd\H os and P.~Tur\'an.  On some sequences of integers.
{\em J. London Math. Soc.} {\bf 11} (1936), 261--264.

\bibitem{FK1}
N.~Frantzikinakis and B.~Kra.
Polynomial averages converge to the product of the integrals.
{\it Isr. J. Math.} 148 (2005), 267-276.

\bibitem{FK2}
N.~Frantzikinakis and B.~Kra.
Convergence of multiple ergodic averages for some
commuting transformations.
{\em Erg. Th. \& Dyn. Sys.} 25 (2005), 799-809.

\bibitem{FK3}
N.~Frantzikinakis and B.~Kra.
Ergodic averages for independent polynomials and applications.
To appear, {\em J. Lond. Math. Soc.}

\bibitem{F3}
H.~Furstenberg.
Strict ergodicity and transformations of the torus.
{\em Amer. J. Math.} {\bf 83}  (1961), 573--601.

\bibitem{F}
H.~Furstenberg.
Ergodic behavior of diagonal measures and a theorem of Szemer\'edi
  on arithmetic progressions.
{\em J. d'Analyse Math.} {\bf 31}  (1977), 204--256.

\bibitem{Fbook}
H.~Furstenberg.
{\em Recurrence in Ergodic Theory and Combinatorial Number Theory}.
Princeton University Press, Princeton, New Jersey, 1981.

\bibitem{F6}
H.~Furstenberg.  Nonconventional ergodic averages.  
{\em Proc. Sympos. Pure Math.} {\bf 50} (1990), 43--56.  

\bibitem{FuK1}
H.~Furstenberg and Y.~Katznelson.
An ergodic Szemer\'edi theorem for commuting transformation.
{\em J. d'Analyse Math.} {\bf 34} (1979), 275--291.

\bibitem{FuK2}
H.~Furstenberg and Y.~Katznelson.
An ergodic Szemer\'edi theorem for IP-systems and combinatorial theory.
{\em J. d'Analyse Math.} {\bf 45} (1985), 117--268.

\bibitem{FKO}
H.~Furstenberg, Y.~Katznelson and D.~Ornstein.  
The ergodic theoretical proof of Szemerédi's theorem.
{\em Bull. Amer. Math. Soc. (N.S.)} {\bf 7} (1982), 527--552.

\bibitem{FW}
H.~Furstenberg and B.~Weiss.
A mean ergodic theorem for $\frac{1}{N} \sum_{n=1}^N
  f({T}^nx)g({T}^{n^2}x)$.
{\em Convergence in Ergodic Theory and Probability}
(Eds.:Bergelson, March, Rosenblatt). Walter de Gruyter \& Co, Berlin, New
  York (1996), 193--227

\bibitem{gowers}
T.~Gowers.
A new proof of Szemer\'edi's theorem.
{\em Geom. Funct. Anal.} {\bf 11} (2001), 465-588; 
Erratum ibid. {\bf 11} (2001), 869.

\bibitem{green}
B.~Green.  
A Szemer\'edi-type regularity lemma in abelian groups.  
{\em Geom. Funct. Anal.} {\bf 15} (2005), 340--376.  
 
\bibitem{GT}
B.~Green and T.~Tao.
The primes contain arbitrarily long arithmetic progressions.
To appear, {\em Annals of Math.}


\bibitem{GT2}
B.~Green and T.~Tao.
An inverse theorem for the Gowers $U^3$ norm.  Preprint, 2005.

\bibitem{GT3}
B.~Green and T.~Tao.
Quadratic uniformity of the M\"obius function.  Preprint, 2005.

\bibitem{GT4}
B.~Green and T.~Tao.
Linear equations in primes.  Preprint, 2006.

\bibitem{H} P. Hall.  A contribution to the theory of groups of prime-power 
order. {\em Proc. London Math. Soc.} (2), \textbf{36} (1933), 29--95.

\bibitem{HVN}
P.~R.~Halmos and J.~von Neumann.  Operator methods in classical 
mechanics, II.  {\em Ann. of Math.}  {\bf 43} (1942), 332--50.

\bibitem{Host}
B.~Host.
Convergence of multiple ergodic averages.
To appear, Proceedings of school on ``Information and Randomness,''
Chili.

\bibitem{HK1}
B.~Host and B.~Kra.
Convergence of Conze-Lesigne averages.
{\em Erg. Th. \& Dyn. Syst.} {\bf 21}  (2001), 493-509.

\bibitem{HK3}
B.~Host and B.~Kra.
Averaging along cubes.
In "Dynamical Systems and
Related Topics,''
Cambridge University Press, 2004, 123--144.

\bibitem{HK4}
B.~Host and B.~Kra.
Nonconventional ergodic
averages and nilmanifolds.  {\em Ann. of Math.} {\bf 161}
(2005), 397--488.

\bibitem{HK5}
B.~Host and B.~Kra. Convergence of polynomial ergodic averages.
{\it Isr. J. Math.} {\bf 149} (2005), 1--19.

\bibitem{johnson}
M.~Johnson.  
Convergence of polynomial ergodic averages of several 
variables for some commuting transformations.  
Preprint, 2006.  

\bibitem{Khintchine}
A.~Y.~Khintchine.
Eine Versch\"{a}rfung des Poincar\'{e}schen
"Wiederkehrsatzes."
{\em Comp. Math.} {\bf 1} (1934), 177--179.

\bibitem{KVN}
B.~O.~Koopman and J.~von Neumann.  Dynamical systems of 
continuous spectra.  {\em Proc. Nat. Acad. Sci. U.S.A.}  {\bf 18} 
(1932), 255-63.  

\bibitem{Kra2}
B.~Kra.
The Green-Tao Theorem on arithmetic progressions in the primes:
an ergodic point of view.  {\em Bull. Amer. Math. Soc.}
{\bf 43} (2006), 3--23.

\bibitem{Kra1}
B.~Kra
From combinatorics to ergodic theory and back again.  
In {\em Proceedings of the International Congress 
of Mathematicians}, Madrid, 2006.  

\bibitem{KN}
L.~Kuipers and N.~Niederreiter.
{\em Uniform distribution of sequences}.  John
Wiley \& Sons, New York, 1974.

\bibitem{La}
M.~Lazard.  
Sur certaines suites d'\'el\'ements dans les groupes 
libres et leurs extensions. 
{\em C. R. Acad. Sci. Paris} {\bf 236} (1953), 36--38.

\bibitem{L}
A. Leibman. Polynomial sequences in groups. {\em J. of Algebra} 
\textbf{201} (1998), 189--206.

\bibitem{Leibman2}
A.~Leibman.
Pointwise convergence of ergodic averages for
polynomial sequences of translations on a nilmanifold.
{\em Erg. Th. \& Dyn. Sys.} (2005), 201-213.

\bibitem{Leibman1}
A.~Leibman.  Convergence of multiple ergodic averages along
polynomials of several variables.   {\em Isr. J. Math.}
{\bf 146} (2005), 303--316.

\bibitem{Lesigne}
E.~Lesigne.  Sur une nil-vari\'et\'e, les parties minimales
associe\`ee \'a une translation sont uniquement ergodiques.
{\em Erg. Th. \& Dyn. Sys.} {\bf 11} (1991), 379--391.

\bibitem{Parry}
W.~Parry.
Ergodic properties of affine transformations and flows on nilmanifolds.
{\em Amer. J. Math.} {\bf 91} (1969), 757--771.


\bibitem{P}
J. Petresco. Sur les commutateurs.  {\em Math. Z.}
 \textbf{61} (1954). 348--356.


\bibitem{poincare}
H.~Poincar\'e. 
Les m\'ethodes nouvelles de la m\'ecanique c\'eleste, I (1892), 
II (1893), and III (1899), Gathiers-Villars, Paris.  


\bibitem{Ratner}
M.~Ratner.  
On Raghunathan's measure conjecture.  {\em Ann. Math.} 
\textbf{134} (1991), 545--607.

\bibitem{sarkozy}
A.~S\'ark\"ozy.
On difference sets of integers I.
{\em Acta Math. Acad. Sci. Hungar.} {\bf 31}
(1978), 125--149.

\bibitem{sarkozy2}
A.~S\'ark\"ozy.
On difference sets of integers III.
{\em Acta Math. Acad. Sci. Hungar.} {\bf 31}
(1978), 355--386.

\bibitem{Shah}
N.~Shah.  
Invariant measures and orbit closures on homogeneous 
spaces for actions of subgroups.  {\em Lie groups 
and ergodic theory} (Mumbai, 1996) Tata Inst. Fund. Res., Bombay 
(1998), 229--271.  

\bibitem{S} E.~Szemer\'edi.
On sets of integers containing no $k$ elements in arithmetic
progression. {\em  Acta Arith.} \textbf{27} (1975), 199--245.

\bibitem{vonneumann}  J.~von Neumann.
Proof of the quasi-ergodic hypothesis.  
{\em Proc. Nat. Acad. Sci. USA} {\bf 18} (1932), 70--82.  

\bibitem{ziegler1}
T.~Ziegler.
A non-conventional ergodic theorem for a nilsystem.
{\em Erg. The. Dyn. Sys.} {\bf 25} (2005), 1357--1370.

\bibitem{ziegler2}
T.~Ziegler.
Universal Characteristic Factors and Furstenberg Averages.  To appear,
{\em J. Amer. Math. Soc.}

\end{thebibliography}

\end{document}